\documentclass[a4paper,12pt]{scrartcl}

\usepackage[utf8]{inputenc}

\usepackage{amsthm,amsmath,amsfonts,amssymb}
\numberwithin{equation}{section}  
\usepackage{mathrsfs}
\usepackage{todonotes}

\usepackage{tikz-cd}

\usepackage{mathtools}	 
\usepackage{bm}	  
\usepackage{bbm}	 

\usepackage{tensor}	 
\usepackage[colorlinks=true]{hyperref} 
\usepackage[all]{hypcap}       
 
\newcommand{\C}{\mathbb{C}}
\newcommand{\R}{\mathbb{R}}
\newcommand{\ii}{\operatorname{i}}

\newcommand{\vol}{\operatorname{vol}}

\newcommand{\Ric}{\operatorname{Ric}}
\newcommand{\del}{\partial}

\newcommand{\PP}{\mathbb{P}}

\newcommand{\cX}{\mathcal{X}}
\newcommand{\Bl}{\operatorname{Bl}}
\newcommand{\rel}{\operatorname{rel}}
\newcommand{\res}{\operatorname{Res}}
\newcommand{\cW}{\mathcal{W}}

\newtheorem{thm}{Theorem}[section]
 
\newtheorem{prop}[thm]{Proposition}

\newtheorem{cor}[thm]{Corollary}

\theoremstyle{definition}
\newtheorem{definition}[thm]{Definition}

\theoremstyle{remark}
\newtheorem{exm}[thm]{Example}
\newtheorem{rmk}[thm]{Remark}

\title{$K$-stability and large \\complex structure limits}
\author{Jacopo Stoppa}
 
\date{\today}

\begin{document}

\maketitle

\begin{abstract} We discuss how, under suitable assumptions, a K\"ahler test configuration admits a mirror Landau-Ginzburg model, giving a corresponding expression for the Donaldson-Futaki invariant as a residue pairing. We study the general behaviour of such mirror formulae under large scaling of the K\"ahler form. We exploit the observation that this scaling trivially preserves $K$-stability, but takes the mirror Landau-Ginzburg model to a large complex structure limit. In certain cases the mirror formulae for the Donaldson-Futaki invariant simplify in this limit. We focus on a special type of limiting behaviour, when the Donaldson-Futaki invariant concentrates at a single critical point of the Landau-Ginzburg potential, and show that this leads to new formulae for the Donaldson-Futaki invariant in terms of theta functions on the mirror. We provide a main application, which shows that such limiting behaviour actually occurs for test configurations in several nontrivial examples, both toric and non-toric, in the case of slope (in)stability for polarised surfaces.
\end{abstract}

\section{Introduction}
Suppose, initially, that $X$ is a nonsingular complex Fano surface (i.e. a del Pezzo), endowed with a maximally degenerate (namely singular, nodal) anticanonical divisor $D$; or $X$ is a toric manifold of arbitrary dimension $n$, with its toric boundary. Note that, in the surface case, $D$ is necessarily the union of irreducible rational curves, $D = \sum^r_{i=1} D_i$, which are smooth for $r > 1$. The fundamental example is a toric Fano manifold endowed with a fixed toric structure and its toric boundary. We also fix a K\"ahler class $[\omega]$ on $X$.

As explained e.g. in \cite{GrossHackingKeel_LCY}, Section 0.5.3 and \cite{KatzKontPant}, Section 1, the mirror to $(X, D, [\omega])$ is an affine manifold $Y_{\omega}$, depending on the K\"ahler class $[\omega]$ for fixed $(X,D)$, together with a nonconstant morphism $W_{\omega}\!: Y \to \C$ and a holomorphic volume form $\Omega_Y$, i.e. a Landau-Ginzburg model,
\begin{equation*}
(X, D, [\omega])\,\big|\, (Y, W, \Omega_Y).
\end{equation*}
The ring of global functions $H^0(Y, \mathcal{O}_Y)$ is endowed with canonical generators $\theta_{D_i}$, corresponding to the boundary components on $X$, and one has
\begin{equation}\label{LGpotentialIntro}
W = \theta_{D_1} + \cdots + \theta_{D_n}.
\end{equation}    
\begin{exm}\label{GiventalFormula} In the toric Fano case, by the classical work of Givental \cite{Givental_toric}, there is a fixed isomorphism $Y \cong (\C^*)^n$ (independent of $[\omega]$), and the theta functions and Landau-Ginzburg potential are given by
\begin{equation*}
\theta_{D_i} = a_i( \omega ) x^{v_i},\,W = \sum^n_{i=1} a_i( \omega ) x^{v_i} 
\end{equation*}
where $v_i$ is the primitive generator of the ray dual to $D_i$ and the $a_i$ are certain coefficients, uniquely determined, up to rescalings of the torus variables, by the condition that, for any integral linear relation
\begin{equation*}
\sum d_i v_i = 0,
\end{equation*}
corresponding to a unique curve class $[C]$ such that $D_i . [C] = d_i$, we have
\begin{equation*}
\prod_i a^{d_i}_i = e^{-2\pi \int_C\omega}. 
\end{equation*}
\end{exm}
In the work \cite{ScarpaStoppa_complexified}, with C. Scarpa, we raised various general questions concerning the behaviour of $K$-stability, the fundamental algebro-geometric stability notion for Fanos and more general polarised manifolds (see e.g. \cite{Donaldson_stabilitySurvey}), with respect to the mirror correspondence between pairs $(X, D)$ and Landau-Ginzburg models. The present paper is motivated by a specific problem.
\newline

\emph{Fix a pair $(X, D)$ with maximally degenerate anticanonical $D \subset X$. Is it possible to find an algebro-geometric characterisation for the mirror Landau-Ginzburg models $W_{\omega}\!: Y_{\omega} \to \C$ corresponding to $K$-(semi, poly)stable K\"ahler classes $[\omega]$ on $(X, D)$?} 
\newline

More precisely, by the discussion above, we are effectively asking for \emph{conditions on the theta functions $\theta_{D_i}$ that cut out the locus of Landau-Ginzburg models $(Y, W)$ mirror to a $K$-(semi, poly)stable $(X, D, [\omega])$ (endowed with an arbitrary K\"ahler class)}. Note that one should allow various notions of $K$-stability such as uniform, log or relative $K$-(semi, poly)stability. The present paper contains some first results towards this general problem. 

In Section \ref{GeneralSec} we discuss how, under suitable assumptions, a K\"ahler test configuration admits a mirror Landau-Ginzburg model, giving a corresponding \emph{expression for the Donaldson-Futaki invariant as a residue pairing computed on the mirror}, see \eqref{DFResPairing} and \eqref{mirrorFutakiIntro}. 

We first study the general behaviour of such mirror formulae under \emph{large scaling of the K\"ahler form}, see Propositions \ref{OscillatoryDFIntroProp} and \ref{OscillatoryResiduesProp}. 

We then exploit the observation that \emph{scaling the K\"ahler form trivially preserves $K$-stability, while this scaling action is highly nontrivial on the mirror, giving one way of approaching a large complex structure limit of the Landau-Ginzburg model. In certain cases the mirror formulae for the Donaldson-Futaki invariant simplify in this limit}.  

We focus on a special type of limiting behaviour, when the Donaldson-Futaki invariant \emph{concentrates at a single critical point of the Landau-Ginzburg potential} (in a sense made precise in \eqref{SpecialLaplaceIntro}), and show that this leads to \emph{new formulae for the Donaldson-Futaki invariant, in terms of theta functions on the mirror, evaluated at the critical point}, see Propositions \ref{MainIdentityIntroProp}, \ref{MainIdentityIntroAlternative} and Corollary \ref{SpecialTestProp}.   

We provide a main application, Theorem \ref{SlopeThm}, which shows that \emph{such limiting behaviour actually occurs for test configurations in several nontrivial examples, both toric and non-toric, in the case of slope (in)stability for polarised surfaces} (a special case of $K$-(in)stability studied by Ross-Thomas \cite{RossThomas_obstruction}). 

As we will recall, in these cases, $(X, [\omega])$ is (potentially) $K$-destabilised by certain special complex submanifolds $Z \subset X$, such that the structure sheaf $\mathcal{O}_Z$ satisfies a ``slope inequality".  From the viewpoint of homological mirror symmetry (see \cite{KatzKontPant}, Section 1), thinking of $\mathcal{O}_Z \in D^b(X)$ as a B-cycle, we may ask if there is a corresponding A-cycle (that is, roughly, a Lagrangian submanifold) on the LG model $Y$ mirror to $X$ with special properties. Indeed, Theorem \ref{SlopeThm} provides a distinguished critical point of the LG potential satisfying a ``mirror slope inequality". On the other hand, it is known that isolated nondegenerate critical points correspond to certain noncompact Lagrangian submanifolds of $Y$, known as Lagrangian thimbles (see e.g. \cite{Auroux_delPezzo}, Section 4.1). This suggests that, in the situation of Theorem \ref{SlopeThm}, there is a ``destabilising A-cycle" corresponding to $\mathcal{O}_Z$ given by a Lagrangian thimble, see Remark \ref{ThimbleRmk}. This observation applies more generally, e.g. Proposition \ref{MainIdentityIntroProp} suggests the existence of Lagrangian thimbles with special numerical properties forced by $K$-(in)stability.\\

\noindent{\textbf{Acknowledgements.}} I would like to thank Mario Garcia-Fernandez, Sohaib Khalid, Riccardo Ontani, Annamaria Ortu, Carlo Scarpa, Lars Martin Sektnan, Nicol\`o Sibilla and especially Giulio Codogni, Ruadha\'i Dervan, Yuji Odaka and Cristiano Spotti for conversations on the topic of this paper. I am very grateful to the anonymous Referees for many important suggestions and corrections.   

\section{General argument}\label{GeneralSec}
 
Suppose $\pi\!: (\mathcal{X}, \omega_{\mathcal{X}}) \to \PP^1$ is a smooth compactified toric K\"ahler test configuration for a compact toric K\"ahler manifold $(X,  \omega)$ of dimension $n$. When $[\omega] = c_1(L)$, this means roughly that $\pi$ restricts to a torus-equivariant polarised flat $\C^*$-degeneration of $(X, L)$ (see \cite{Legendre_localised} for more details in the general case). 
\begin{rmk} The toric assumption is convenient for the present discussion, especially since we can rely on the very general mirror symmetry results of Coates, Corti, Iritani and Tseng \cite{CoatesCortiIritani_hodge}, but ideally one should be able to work with the theta functions constructed in the Gross-Siebert programme \cite{GrossHackingSiebert_theta}. Indeed we do that in this paper for certain non-toric del Pezzos, as we explain in Section \ref{GeneralExmSec}. 
\end{rmk}
To each such test configuration one associates a Donaldson-Futaki weight, or invariant, denoted by $\operatorname{DF}(\mathcal{X}, \omega_{\cX})$, and $K$-semistability of $(X, \omega)$ is defined as the condition $\operatorname{DF}(\mathcal{X}, \omega_{\cX}) \geq 0$ for all $(\mathcal{X}, \omega_{\cX})$ (this coincides with the more usual definition of $K$-semistability involving normal test configuratons, in the toric case, by \cite{DervanRoss_Kstab}). 

By the results of Odaka and Wang \cite{Odaka_blowup, Wang_GIT}, the Donaldson-Futaki invariant can be expressed as the Poincar\'e pairing
\begin{equation*}
\operatorname{DF}(\mathcal{X}, \omega_{\cX}) = \int_{\mathcal{X}} (\mathcal{A}_{\omega_{\mathcal{X}}})^n \cup \mathcal{B}_{\omega_{\mathcal{X}}}
\end{equation*}
where the differential forms $\mathcal{A}_{\omega_{\mathcal{X}}}$, $\mathcal{B}_{\omega_{\mathcal{X}}}$ are given by 
\begin{align*}
\mathcal{A}_{\omega_{\mathcal{X}}} = \omega_{\mathcal{X}},\,\mathcal{B}_{\omega_{\mathcal{X}}} = \frac{n \mu}{n+1} \omega_{\mathcal{X}} -  \Ric(\omega_{\mathcal{X}}) +  \pi^*\omega_{FS},
\end{align*}
and we set
\begin{equation*}
\mu = -\frac{K_X . [\omega]^{n-1}}{[\omega]^n}. 
\end{equation*}

Turning to mirror-symmetric constructions, we recall that, since $\cX$ is toric, the general results of \cite{CoatesCortiIritani_hodge} provide a mirror family to $\cX$ as a torus fibration $\mathcal{Y} \to \mathcal{M}$, with fibres isomorphic to $(\C^*)^{n+1}$, endowed with a regular function $\cW\!: \mathcal{Y} \to \C$, the Landau-Ginzburg (LG) potential. The choice of the K\"ahler class $[\omega_{\cX}]$ determines a point on the base of the fibration, which is denoted by $\tau^{-1}(0)$ in the notation of \cite{CoatesCortiIritani_hodge}, Theorem 4.28 (where the variable $Q$ appearing there satisfies $Q^d = \exp(-2 \pi \int_d \omega_{\mathcal{X}}$), and so a regular function 
\begin{equation*}
\cW(\omega_{\cX}):= \cW|_{\mathcal{Y}_{\tau^{-1}(0)}}\!:\mathcal{Y}_{\tau^{-1}(0)} \to \C.
\end{equation*} 
\begin{rmk}
The construction of the fibration $\mathcal{Y} \to \mathcal{M}$ in \cite{CoatesCortiIritani_hodge} depends on the choice of a finite set $G$ of ``unfolding parameters". This is required when $\cX$ is not weak Fano: then we need to choose an appropriate $G$ so that the mirror map $\tau$ is a local isomorphism (or at least, such that $\tau$ is surjective and all points of the fibre $\tau^{-1}(0)$ are related by coordinate changes of the torus). We fix such a choice of $G$ in the following. Also, we need to choose the K\"ahler class $[\omega_{\cX}]$ to be sufficiently large so that the mirror map $\tau$ is convergent (according to the results of \cite{CoatesCortiIritani_hodge}, Section 7). This is not restrictive for our applications since we will work in the limit $k[\omega_{\cX}]$ for $k \to + \infty$.

Moreover, there is a canonical extension to a singular fibration $\widehat{\mathcal{Y}} \to \widehat{\mathcal{M}}$ for a partial compactification $\mathcal{M} \subset \widehat{\mathcal{M}}$, with fibres degenerating over the so-called large complex structure limit point in $\widehat{\mathcal{M}}$. 

It is also important to point out that the results of \cite{CoatesCortiIritani_hodge} hold in the stronger equivariant setting, with respect to the action of the full torus on $\cX$. Here we consider the non-equivariant limit, but see Remark \ref{EquivRmk}. 
\end{rmk}

Let us describe the basic properties satisfied by the mirror construction for our purposes. Write $\operatorname{Jac}(\cW(\omega_{\cX})) = \Gamma(\mathcal{O}_{\operatorname{Crit}(\cW(\omega_{\cX}))})$ for the Jacobi ring of the regular function $\cW(\omega_{\cX})$. Then, there exists a canonically defined rank $1$ free module over $\operatorname{Jac}(\cW(\omega_{\cX}))$, denoted by $\operatorname{GM}(\cW(\omega_{\cX}))|_{z=0}$ (a quotient of the module of relative $n+1$-forms), endowed with a natural bilinear form $P_{\cW(\omega_{\cX})}|_{z=0}$, called the \emph{residue pairing}, and a ``mirror map" isomorphism of vector spaces
\begin{equation*}
\Theta_{\omega_{\cX}}\!: \operatorname{GM}(\cW(\omega_{\cX}))|_{z = 0} := H^{n+1}\left(\Omega^{\bullet}_{\widehat{ \mathcal{Y}} /\widehat{ \mathcal{M}}}, d \cW(\omega_{\cX})\wedge\right) \to H^*(\mathcal{X}, \C),
\end{equation*}
which intertwines the Poincar\'e pairing and the residue pairing, see \cite{CoatesCortiIritani_hodge}, Section 6.

Thus, when $(X, [\omega])$ is $K$-semistable, we find that for all the LG potentials $\mathcal{W}(\omega_{\mathcal{X}})$ obtained from $(\cX, \omega_{\cX})$ as above we must have
\begin{equation}\label{DFResPairing}
P_{\cW(\omega_{\cX})}|_{z = 0}(\Theta^{-1}_{\omega_{\cX}}([\mathcal{A}_{\omega_{\mathcal{X}}}]^n), \Theta^{-1}_{\omega_{\cX}}([\mathcal{B}_{\omega_{\mathcal{X}}}])) = \operatorname{DF}(\mathcal{X}, \omega_{\cX})\geq 0.
\end{equation}
This is a nontrivial constraint on a Landau-Ginzburg potential $ W(\omega )$ mirror to a $K$-semistable toric manifold $(X, [\omega])$. We will discuss the interpretation of this constraint in Remark \ref{LGtestRmk}.

Next, note that a choice of a generator induces an algebra structure from $\operatorname{Jac}(\cW(\omega_{\cX}))$ to $\operatorname{GM}(\cW(\omega_{\cX}))|_{z = 0}$. If we choose the generator $\Theta^{-1}_{\omega_{\cX}}(1)$ (where $1 \in H^*(\mathcal{X}, \C)$), then $\Theta_{\omega_{\cX}}$ is an isomorphism between $\operatorname{GM}(\cW(\omega_{\cX}))|_{z = 0}\cong\operatorname{Jac}(\cW(\omega_{\cX})$ and the \emph{small quantum cohomology ring} $(H^*(\mathcal{X}, \C), *_{Q(\omega_{\cX})})$ (i.e. for the quantum parameter $Q = Q(\omega_{\cX})$ such that $Q^d = \exp(-2 \pi \int_d \omega_{\mathcal{X}})$, see \cite{CoatesCortiIritani_hodge}, Section 4.) 

This fundamental fact fact allows to study the behaviour of the ``stability condition" \eqref{DFResPairing} under large scalings of the K\"ahler form. That is, we note that $(X, [k\omega])$ is trivially $K$-semistable for all $k > 0$, but the change $\mathcal{W}(\omega_{\mathcal{X}})\mapsto\mathcal{W}(k\omega_{\mathcal{X}})$ is highly nontrivial. Indeed, on the Landau-Ginzburg side, $k \to \infty$ corresponds to a large complex structure limit. Note that we have
\begin{align*}
& [\mathcal{A}_{k\omega_{\mathcal{X}}} ]^n = [\mathcal{A}_{k\omega_{\mathcal{X}}}] \cup \cdots \cup [\mathcal{A}_{k\omega_{\mathcal{X}}}] = \lim_{k \to \infty}\,[\mathcal{A}_{k\omega_{\mathcal{X}}}] *_{Q(k\omega_{\cX})} \cdots *_{Q(k\omega_{\cX})} [\mathcal{A}_{k\omega_{\mathcal{X}}}],
\end{align*}
since the small quantum cohomology product $*_{Q(k\omega_{\cX})}$ converges to the usual cup product for $k \to \infty$. Moreover, the product of functions in the algebra $\operatorname{GM}(\mathcal{W}(k\omega_{\mathcal{X}}))|_{z = 0} \cong \operatorname{Jac}(\cW(k\omega_{\cX}))$ is intertwined with $*_{Q(k\omega_{\cX})}$ by the mirror isomorphism $\Theta_{k\omega_{\cX}}$, so we find
\begin{align*}
&\big(k^{-1}\Theta^{-1}_{k\omega_{\cX}}(\mathcal{A}_{k\omega_{\mathcal{X}}})\big)^n = k^{-n}\Theta^{-1}_{k\omega_{\cX}}\big([\mathcal{A}_{k\omega_{\mathcal{X}}}] *_{Q(k\omega_{\cX})} \cdots *_{Q(k\omega_{\cX})} [\mathcal{A}_{k\omega_{\mathcal{X}}}]\big)\\
& = k^{-n}\Theta^{-1}_{k\omega_{\cX}}\big( [\mathcal{A}_{k\omega_{\mathcal{X}}} ]^n \big) + O(k^{-n-1}).
\end{align*}
We can now use the scaling behaviour $\operatorname{DF}(\mathcal{X}, [k\omega_{\mathcal{X}}]) = k^n \operatorname{DF}(\mathcal{X}, [\omega_{\mathcal{X}}])$, together with the fact that the form $\mathcal{B}_{\omega_{\mathcal{X}}}$ is scale invariant, to obtain the following.
\begin{prop}\label{OscillatoryDFIntroProp} 
For $k >0$ we have an expansion
\begin{align*} 
\operatorname{DF}(\mathcal{X}, [\omega_{\mathcal{X}}]) = P_{\cW(k\omega_{\cX})}|_{z = 0}(\big(k^{-1}\Theta^{-1}_{k\omega_{\cX}}(\mathcal{A}_{k\omega_{\mathcal{X}}})\big)^n, \Theta^{-1}_{k\omega_{\cX}}([\mathcal{B}_{k\omega_{\mathcal{X}}}])) + O(k^{-1}).
\end{align*}
In particular, if $(\mathcal{X}, [\omega_{\mathcal{X}}])$ is strictly destabilising, we must have
\begin{equation*}
P_{\cW(k\omega_{\cX})}|_{z = 0}(\big(k^{-1}\Theta^{-1}_{k\omega_{\cX}}(\mathcal{A}_{k\omega_{\mathcal{X}}})\big)^n, \Theta^{-1}_{k\omega_{\cX}}([\mathcal{B}_{k\omega_{\mathcal{X}}}])) < 0
\end{equation*}
for all sufficiently large $k$.
\end{prop}
In order to make this observation more explicit we can use the basic fact that the mirror fibration $\mathcal{Y} \to \mathcal{M}$ is in fact trivial, although it is not canonically trivialised (see \cite{CoatesCortiIritani_hodge}, Section 4). Fixing a choice of trivialisation, we can regard the LG potential $\cW(k\omega_{\cX})$ as a regular function defined on the algebraic torus $(\C^*)^{n+1}$. Similarly, the algebra $\operatorname{GM}(\cW(k\omega_{\cX}))|_{z = 0}\cong\operatorname{Jac}(\cW(k\omega_{\cX}))$ (with isomorphism induced by $\Theta^{-1}_{k\omega_{\cX}}(1)$) can be regarded as a quotient of the algebra of regular functions on $(\C^*)^{n+1}$.

\begin{definition} Fix a choice of trivialisation of $\mathcal{Y} \to \mathcal{M}$ and a K\"ahler form $\omega_{\cX}$. We denote by $\vartheta_{\omega_{\cX},k}$ (a ``generalised theta function") any lift of $k^{-1}\Theta^{-1}_{k\omega_{\cX}}([k\omega_{\cX}])$ to a regular function on $(\C^*)^{n+1}$. Similarly, we write $\mathcal{W}_{\operatorname{rel},k} = \mathcal{W}_{\operatorname{rel},k}(k\omega_{\cX})$ (a ``generalised relative LG potential") for any lift of $\Theta^{-1}_{k\omega_{\cX}}(c_1(K^{\vee}_{\mathcal{X}/\PP^1}))$ to a regular function on $(\C^*)^{n+1}$. Finally we set $\cW_k = \cW(k\omega_{\cX})$.
\end{definition}  
\begin{rmk} This notation can be partially justified by the fact that, at least when $\cX$ is Fano, there is a canonical choice of the lift $\mathcal{W}_{\operatorname{rel},k}(k\omega_{\cX})$ and a set of distinguished choices for $\vartheta_{\omega_{\cX},k}$ (for a fixed trivialisation of $\mathcal{Y} \to \mathcal{M}$). Indeed, by the results of \cite{CoatesCortiIritani_hodge}, if $\cX$ is \emph{weak} Fano, with toric boundary $\mathcal{D} = \mathcal{D}_1 + \cdots + \mathcal{D}_{m}$, then the LG potential can be written in the form
\begin{equation*}
\cW_k = \mathcal{W}(k\omega_{\cX}) = \sum^{m}_{i=1} \vartheta_{\mathcal{D}_i,k}(k\omega_{\cX})
\end{equation*}
for certain unique ``theta functions" $\vartheta_{\mathcal{D}_i,k}(k\omega_{\cX})$ corresponding to the boundary components $\mathcal{D}_i$ (depending on the K\"ahler class $[k\omega_{\mathcal{X}}]$) such that, following the notation of Example \ref{GiventalFormula}, we have
\begin{equation*}
\vartheta_{\mathcal{D}_i,k}(k\omega_{\cX}) = \theta_{\mathcal{D}_i}(k\omega_{\cX})(1+ O(k^{-1})) 
\end{equation*}
(note that the quantity $O(k^{-1})$ does not depend on the torus variables $x_j$, ensuring that $\vartheta_{\mathcal{D}_i,k}$ remains a monomial in the variables $x_j$). In general, if $\mathcal{D} = \sum_{i} b_i \mathcal{D}_i$, we set $\vartheta_{\mathcal{D}} = \sum_{i} b_i \vartheta_{\mathcal{D}_i}$. Then, if $\cX$ is actually Fano, writing 
\begin{equation*}
[\omega_{\mathcal{X}}] = \sum^{m}_{i=1} w_i [\mathcal{D}_i],\,w_i \in \R, 
\end{equation*}
we can choose
\begin{align*}
&\vartheta_{\omega_{\mathcal{X}},k} = \sum^{m}_{i=1} w_i \vartheta_{\mathcal{D}_i, k},\,\mathcal{W}_{\operatorname{rel},k} = \sum^{m}_{i=1} \vartheta_{\mathcal{D}_i,k} - \vartheta_{\pi^*\mathcal{O}_{\PP^1}((0) + (\infty)),k} = \mathcal{W}(k\omega_{\cX}) -\vartheta_{\mathcal{X}_0,k}-\vartheta_{\mathcal{X}_{\infty},k}, 
\end{align*}
where $\mathcal{X}_0$, $\mathcal{X}_{\infty}$ are the central fibre and the divisor at infinity respectively.
\end{rmk}
\begin{rmk}
Similarly, working on a Fano $X$ with toric boundary $D = D_1 + \cdots + D_n$, we write
\begin{align*}
W_k = W(k[\omega_X]) = \sum^n_{i=1} \vartheta_{D_i, k},\,[\omega_X] = \sum^n_{i = 1} a_i [D_i],\,\vartheta_{\omega_X, k} = \sum^n_{i=1} a_i \vartheta_{D_i, k}.
\end{align*}
\end{rmk}

With this notation in place, we can use the explicit description of the residue pairing for regular functions on $(\C^*)^{n+1}$ (see \cite{CoatesCortiIritani_hodge}, Section 6) together with Proposition \ref{OscillatoryDFIntroProp} to obtain 
\begin{align*}  
&\operatorname{DF}(\mathcal{X}, [\omega_{\mathcal{X}}]) = P_{\cW(k\omega_{\cX})}|_{z = 0}(\big(k^{-1}\Theta^{-1}_{k\omega_{\cX}}(\mathcal{A}_{k\omega_{\mathcal{X}}})\big)^n, \Theta^{-1}_{k\omega_{\cX}}([\mathcal{B}_{k\omega_{\mathcal{X}}}])) + O(k^{-1})\\
&\nonumber= \res_{\mathcal{W}(k\omega_{\mathcal{X}})}(\vartheta^n_{\omega_{\mathcal{X}},k}, \frac{n \mu}{n+1}\vartheta_{\omega_{\mathcal{X}},k} - \mathcal{W}_{\operatorname{rel},k}) + O(k^{-1})
\end{align*}
where the classical residue pairing (as defined e.g. in \cite{GriffithsHarris}, p. 659 and \cite{Hertling}, Section 10.4) is given by
\begin{align*}  
&\res_{\mathcal{W}(k\omega_{\mathcal{X}})}(\vartheta^n_{\omega_{\mathcal{X}},k}, \frac{n \mu}{n+1}\vartheta_{\omega_{\mathcal{X}},k} - \mathcal{W}_{\operatorname{rel},k}) \\
&= \left(\frac{1}{2\pi \ii}\right)^{n+1}\int_{\Gamma} \frac{ \vartheta^n_{\omega_{\mathcal{X}},k} (\frac{n \mu}{n+1}\vartheta_{\omega_{\mathcal{X}},k} - \mathcal{W}_{\operatorname{rel},k})}{ \prod_{i} x_i\del_{x_i} \mathcal{W}_k } \Omega(k \omega_{\cX}), 
\end{align*}
$\Omega(k \omega_{\cX})$ is a lift of $\Theta^{-1}_{k\omega_{\cX}}(1)$, and the integration cycle is given by
\begin{equation*}
\Gamma = \{|\del_{x_i} \mathcal{W}| = \varepsilon\},\,0 < \varepsilon \ll 1.
\end{equation*}
(with positive orientation). Using the property
\begin{equation*}
\Omega(k \omega_{\cX}) = \Omega + O(k^{-1}),
\end{equation*}
where the standard holomorphic volume form is 
\begin{equation*}
\Omega = \frac{dx_1}{x_1} \wedge \cdots \wedge \frac{d x_{n+1}}{x_{n+1}},
\end{equation*}
shows the following.
\begin{prop} We have
\begin{align}\label{OscillatoryDFIntro} 
\operatorname{DF}(\mathcal{X}, [\omega_{\mathcal{X}}]) = \lim_{k \to \infty} \left(\frac{1}{2\pi \ii}\right)^{n+1}\int_{\Gamma} \frac{ \vartheta^n_{\omega_{\mathcal{X}},k} (\frac{n \mu}{n+1}\vartheta_{\omega_{\mathcal{X}},k} - \mathcal{W}_{\operatorname{rel},k})}{ \prod_{i} x_i\del_{x_i} \mathcal{W}_k } \Omega.
\end{align}
\end{prop}
\begin{rmk}
When the potentials $\mathcal{W}_k$ have nondegenerate critical points, the residue pairing admits a simple expression which will be recalled in Proposition \ref{OscillatoryResiduesProp} below. 
\end{rmk}

\begin{exm}[Degeneration to the normal cone of a point]\label{P1Exm} It seems helpful to work out in detail the basic example obtained when $(X, [\omega])$ is given by $(\PP^1, [\omega_{FS}])$ and the total space of $\cX$ is given by degeneration to the normal cone of a point in $\PP^1$, i.e. $\cX = \Bl_{(p,0)} (\PP^1 \times \PP^1)$, where the second copy of $\PP^1$ is thought of as $\C \cup \{0, \infty\}$. We realise $\cX$ as a toric surface with fan
\begin{equation*}
v_1 = (1, 0),\,v_2 = (1,1),\,v_3 = (0,1),\,v_4=(-1,0),\,v_5=(0,-1).
\end{equation*} 
Dually to the primitive fan generators we have divisors
\begin{equation*}
L_1,\,E,\,L_2,\,H_2,\,H_1,
\end{equation*}
where $H_1 = \{q\} \times \PP^1$, $H_2 = \PP^1 \times \{\infty\}$ for $q \neq p$, $L_1$, $L_2$ are the proper transforms of the corresponding fibres through $(p,0)$ and $E$ is the exceptional divisor. 

The K\"ahler class on the test configuration $\cX$ is given by
\begin{equation*}
[\omega_{\cX}] = H_1 + H_2 - r E,\,r\in(0, 1),
\end{equation*}
and we have, as toric divisors,
\begin{equation*}
\cX_{0} = L_1 + E,\, \cX_{\infty} = H_2.
\end{equation*}
Note that $\cX$ is Fano, i.e. a del Pezzo surface. Thus, we can write the Landau-Ginzburg potential as
\begin{align*}
&\cW_k = \cW(k\omega_{\cX}) = \vartheta_{H_1, k} + \vartheta_{H_2, k} + \vartheta_{L_2, k} + \vartheta_{E, k} + \vartheta_{L_1, k} \\
&= \theta_{H_1, k} + \theta_{H_2, k} + \theta_{L_2, k} + \theta_{E, k} + \theta_{L_1, k} 
\end{align*}
where the theta functions are given by 
\begin{align*}
\theta_{H_1, k} = \frac{e^{-2\pi k}}{y},\,\theta_{H_2, k} =\frac{e^{-2\pi k}}{x},\,\theta_{L_2, k} = y,\,\theta_{E, k} = e^{2\pi k r} x y,\,\theta_{L_1, k} = x, 
\end{align*}
$(x, y)$ denoting coordinates on $(\C^*)^2$. Similarly we have
\begin{align*}
& \cW_{\rel, k} = \cW_k - \vartheta_{\mathcal{X}_0,k}-\vartheta_{\mathcal{X}_{\infty},k} =  \theta_{H_1, k} + \theta_{L_2, k},\\
&\vartheta_{\omega_{\cX}, k} = \theta_{H_1,k} + \theta_{H_2,k} - r \theta_{E,k}. 
\end{align*}
\end{exm}
\begin{rmk}[Test configurations for LG potentials]\label{LGtestRmk} To some extent, $\mathcal{W}_k\!: (\C^*)^{n+1} \to \C$ may be considered as a ``test configuration" for $W_k\!: (\C^*)^n \to \C$. It is easier to see this in the special case when $\cX$ the non-compact toric test-configuration $\mathcal{X}^o$ (such that $\mathcal{X}$ is the canonical compactification of $\mathcal{X}^o$) is induced by an iterated toric blowup of the trivial test-configuration $X \times \C$, where the extra $\C^*$ acts on $\C$ in the standard way (as in the work of Odaka \cite{Odaka_blowup}). Then, $\mathcal{X}^o$ is determined by a refinement of the toric fan of $X \times \C$, and so, as $\cX^o$ is toric, projective over an affine, and contains a torus fixed point, by the results of \cite{CoatesCortiIritani_hodge}, the Landau-Ginzburg potential $\mathcal{W}^o_k$ corresponding to $(\mathcal{X}^o, k \omega_{\cX}|_{\mathcal{X}^o})$ is of the form 
\begin{equation*}
\mathcal{W}^o_k(x, y) = \big(\widetilde{W}_k(x) + y\,\mathcal{U}_k(x, y)\big)\big(1 + O(k^{-1})\big)\!:(\C^*)^n_{x} \times \C^*_y \to \C,\, x = (x_1, \ldots, x_n).
\end{equation*}
Here, $\widetilde{W}_k(x)$ contains the same monomials as $W_k(x)$, and specialises to $W_k(x)$ in the limit $[\omega_{\mathcal{X}}] \to p^*_X[\omega_X]$, while $\mathcal{U}_k(x, y)$ is a \emph{polynomial}, rather than a Laurent polynomial, in $y$ (since the exceptional locus of $\mathcal{X}^o \to X \times \C$ maps to $0 \in \C$), containing \emph{arbitrary} monomials $x^i_l y^j$ with $j \geq 0$, and with coefficients determined by $[k\omega_{\mathcal{X}}]$. So $\mathcal{W}^o_k(x, y)$ specialises to $W_k(x)$ for $y = 0$ and $[\omega_{\mathcal{X}}] \to p^*_X[\omega_X]$, and the Landau-Ginzburg potential of the compactification $\mathcal{X}$ is given by
\begin{equation}\label{WBlowup}
\mathcal{W}_k = \big(\frac{a_k}{y} + \widetilde{W}_k(x) + y\,\mathcal{U}_k(x, y)\big)\big(1 + O(k^{-1})\big)
\end{equation}
for some $a_k \in \C^*$ determined by $[k\omega_{\mathcal{X}}]$. 

Finally, we can think that the ``central fibre" is given by restricting to the hypersurface 
\begin{equation*}
H = \{\del_y \mathcal{W}_k = 0\} \subset (\C^*)^n_{x} \times \C^*_y.
\end{equation*} 
Indeed, according to \eqref{OscillatoryDFIntro}, the Donaldson-Futaki invariant can be computed in an arbitrarily small open neighbourhood of $H$.

In the situation of Example \ref{P1Exm} we compute  
\begin{align*}
&a_k  e^{-2\pi k},\,\widetilde{W}_k(x) = \frac{e^{-2\pi k}}{x} + x,\,\mathcal{U}_k(x, y) = 1 + e^{2\pi k r} x,
\end{align*} 
and the ``central fibre" is given by restriction to $H = \{\mathcal{U}_k\,y^2 = e^{-2\pi k}\}$. 
\end{rmk}

\begin{rmk}[Equivariance]\label{EquivRmk} As we mentioned, the results of \cite{CoatesCortiIritani_hodge} hold equivariantly for the full torus action on $\mathcal{X}$, while in our discussion above we only considered the non-equivariant limit. However, it is possible to preserve some equivariance, namely with respect to the distinguished action of $S^1$ on $\mathcal{X}$ which is part of the structure of a test configuration. This is because, as discussed by Legendre \cite{Legendre_localised}, the forms $\mathcal{A}_{\omega_{\mathcal{X}}}$, $\mathcal{B}_{\omega_{\mathcal{X}}}$ can be lifted to $S^1$-equivariant forms given explicitly by     
\begin{align*}
&\widehat{\mathcal{A}}_{\omega_{\mathcal{X}}} = \omega_{\mathcal{X}} - \langle m , v \rangle,\\
&\widehat{\mathcal{B}}_{\omega_{\mathcal{X}}} = \frac{n \mu}{n+1} (\omega_{\mathcal{X}} - \langle m, v\rangle ) -  (\Ric(\omega_{\mathcal{X}}) - \frac{1}{2}\Delta_{\omega_{\mathcal{X}}} \langle m, v\rangle ) +  (\pi^*\omega_{FS} - \pi^*m_{FS} ),
\end{align*}
where $\omega_{\cX}$ is $S^1$-invariant, $m$, $m_{FS}$ denote equivariant moment maps for the $S^1$-actions on $(\mathcal{X}, \omega_{\mathcal{X}})$, respectively $\PP^1$ endowed with the Fubini-Study form, and $v$ is the infinitesimal generator. The forms $\widehat{\mathcal{A}}_{\omega_{\mathcal{X}}}$, $\widehat{\mathcal{B}}_{\omega_{\mathcal{X}}}$ are equivariantly closed for the $S^1$-action on $\mathcal{X}$ and so define elements of equivariant cohomology, 
\begin{equation*}
[\widehat{\mathcal{A}}_{\omega_{\mathcal{X}}}], \,[\widehat{\mathcal{B}}_{\omega_{\mathcal{X}}}] \in H^*_{\mathbb{T}}(\mathcal{X}). 
\end{equation*}
According to \cite{Legendre_localised} (building on \cite{Odaka_blowup, Wang_GIT}), the Donaldson-Futaki invariant can be expressed as a Poincar\'e pairing of equivariant classes by
\begin{equation*}
\operatorname{DF}(\mathcal{X}, [\omega_{\mathcal{X}}]) = \int_{\mathcal{X}} (\widehat{\mathcal{A}}_{\omega_{\mathcal{X}}})^n \cup \widehat{\mathcal{B}}_{\omega_{\mathcal{X}}}. 
\end{equation*}
We may then repeat our discussion above with respect to a partial specialisation of the equivariant mirror isomorphism of \cite{CoatesCortiIritani_hodge}, given by
\begin{equation*}
\Theta^{S^1}_{\omega_{\cX}}\!: \operatorname{GM}(\cW_{S^1}(\omega_{\cX}))|_{z = 0} \to H^*_{S^1}(\mathcal{X}).
\end{equation*}

\end{rmk}
\subsection{Mirror Donaldson-Futaki invariant}
Let $x_1, \ldots, x_{n+1}$ denote torus coordinates on $(\C^*)^{n+1}$. The following expression follows immediately from \eqref{OscillatoryDFIntro} and the characterisation of the residue pairing in \cite{CoatesCortiIritani_hodge}, Section 6.
 \begin{prop}\label{OscillatoryResiduesProp} If $\mathcal{W}_k = \mathcal{W}(k\omega_{\mathcal{X}})$ has nondegenerate critical points for $k \gg 1$, we have
 \begin{align}\label{OscillatoryResidues} 
\operatorname{DF}(\mathcal{X}, [\omega_{\mathcal{X}}]) &= \lim_{k\to\infty}\sum_{p \in \operatorname{Crit}(\mathcal{W}_k)} \frac{ \vartheta^n_{\omega_{\mathcal{X}},k} \big(\frac{n \mu}{n+1}\vartheta_{\omega_{\mathcal{X}},k} - \mathcal{W}_{\operatorname{rel},k}\big)}{\prod_i x^2_i \det \nabla^2 \mathcal{W}_k}\big|_p. 
\end{align} 
 \end{prop}

Note that, by the same arguments as above using the mirror theorem, we also have
\begin{equation*}
\vol(\mathcal{X}) = \lim_{k\to\infty}\sum_{p \in \operatorname{Crit}(\mathcal{W}_k)} \vol_p(\mathcal{X}),
\end{equation*}
where we set
\begin{equation*}
\vol_p(\mathcal{X}) := \frac{ \vartheta^{n+1}_{\omega_{\mathcal{X}},k}}{\prod_i x^2_i \det \nabla^2 \mathcal{W}_k}\big|_p.
\end{equation*}
As a consequence, we have
\begin{cor} If $\mathcal{W}_k = \mathcal{W}(k\omega_{\mathcal{X}})$ has nondegenerate critical points for $k \gg 1$, we have
 \begin{align}\label{OscillatoryResiduesVol} 
\operatorname{DF}(\mathcal{X}, [\omega_{\mathcal{X}}]) &= \lim_{k\to\infty}\sum_{p \in \operatorname{Crit}(\mathcal{W}_k)}  \vol_{p}(\mathcal{X}) \left(\frac{n\mu}{n+1} - \frac{ \mathcal{W}_{\operatorname{rel},k} }{\vartheta_{\omega_{\mathcal{X}},k} }\right)\big|_{p}. 
\end{align}
 \end{cor}
\begin{exm}\label{P1ExmCont} We continue with the basic Example \ref{P1Exm}. Let us make the special choice $r = \frac{1}{3}$. We have
\begin{equation*}
\nabla \cW_k = \left(-\frac{e^{-2 \pi  k}}{x^2}+e^{\frac{2 \pi  k}{3}} y+1,e^{\frac{2 \pi  k}{3}} x-\frac{e^{-2 \pi  k}}{y^2}+1\right)
\end{equation*}
and a little computation shows that the critical points of $\cW_k$ are cut out by the equations
\begin{align*}
&\left(e^{2 \pi  k} x^2-e^{\frac{2 \pi  k}{3}} x-1\right) \left(e^{\frac{8 \pi  k}{3}} x^3+e^{2 \pi  k} x^2-1\right) = 0,\,y = e^{-\frac{8 \pi  k}{3}} \left(\frac{1}{x^2}-e^{2 \pi  k}\right).
\end{align*}
Let us check \eqref{OscillatoryResidues} in this case. There are five critical points, and using the equations above we can work out their possible asymptotic behaviours
\begin{align*}
& p_1 \sim (-e^{-\frac{2}{3}\pi k} , -e^{-\frac{2}{3}\pi k}),\,p_i \sim (\pm e^{-\pi k}, e^{-\pi k}),\,i = 2, 3,\, p_i \sim (\pm e^{-\pi k}, - e^{-\pi k}),\,i = 4, 5.
\end{align*}
We also have
\begin{align*}
&\vartheta_{\omega_{\cX}, k} = -\frac{1}{3} e^{\frac{2 \pi  k}{3}} x y+\frac{e^{-2 \pi  k}}{x}+\frac{e^{-2 \pi  k}}{y},\\
&\frac{1}{(x y)^2 \det \nabla^2 \mathcal{W}_k} = \frac{e^{4 \pi  k} x y}{4-e^{\frac{16 \pi  k}{3}} x^3 y^3}  
\end{align*}
and so, using $n = 1, \mu = 2$, 
\begin{equation*}
\frac{ \vartheta_{\omega_{\mathcal{X}},k} \big(\vartheta_{\omega_{\mathcal{X}},k} - \mathcal{W}_{\operatorname{rel},k}\big)}{\prod_i x^2_i \det \nabla^2 \mathcal{W}_k} = \frac{\left(e^{\frac{8 \pi  k}{3}} x^2 y+3 e^{2 \pi  k} x y-3\right) \left(e^{\frac{8 \pi  k}{3}} x^2 y^2-3 x-3 y\right)}{9 x \left(4-e^{\frac{16 \pi  k}{3}} x^3
   y^3\right)}.
\end{equation*}
Plugging in our asymptotics, we obtain
\begin{align*}
& \frac{ \vartheta_{\omega_{\mathcal{X}},k} \big(\vartheta_{\omega_{\mathcal{X}},k} - \mathcal{W}_{\operatorname{rel},k}\big)}{\prod_i x^2_i \det \nabla^2 \mathcal{W}_k}\big|_{p_1} \sim \frac{\left(e^{\frac{2 \pi  k}{3}}+6\right) \left(2 e^{\frac{2 \pi  k}{3}}-3\right)}{9 \left(e^{\frac{4 \pi  k}{3}}-4\right)} \to \frac{2}{9},\\
& \frac{ \vartheta_{\omega_{\mathcal{X}},k} \big(\vartheta_{\omega_{\mathcal{X}},k} - \mathcal{W}_{\operatorname{rel},k}\big)}{\prod_i x^2_i \det \nabla^2 \mathcal{W}_k}\big|_{p_i} \sim \pm\frac{1-6 e^{\frac{\pi  k}{3}}}{36 e^{\frac{2 \pi  k}{3}} \mp 9} \to 0,\,i = 2, 3,\\
& \frac{ \vartheta_{\omega_{\mathcal{X}},k} \big(\vartheta_{\omega_{\mathcal{X}},k} - \mathcal{W}_{\operatorname{rel},k}\big)}{\prod_i x^2_i \det \nabla^2 \mathcal{W}_k}\big|_{p_i} \sim \mp\frac{6 e^{\frac{\pi  k}{3}}+1}{36 e^{\frac{2 \pi  k}{3}}\pm9} \to 0,\,i = 4, 5, 
\end{align*}
so we find
\begin{equation*}
\operatorname{DF}(\mathcal{X}, [\omega_{\mathcal{X}}]) =\lim_{k\to\infty} \frac{ \vartheta_{\omega_{\mathcal{X}},k} \big(\vartheta_{\omega_{\mathcal{X}},k} - \mathcal{W}_{\operatorname{rel},k}\big)}{\prod_i x^2_i \det \nabla^2 \mathcal{W}_k}\big|_{p_1} = \frac{2}{9}. 
\end{equation*}
This agrees with the intersection-theoretic formula
\begin{align*}
&\operatorname{DF}(\mathcal{X}, [\omega_{\mathcal{X}}]) = [\omega_{\cX}]^2 + K_{\cX/\PP^1}. [\omega_{\cX}] \\
&= (H_1 + H_2 - \frac{1}{3} E)^2 - (H_1 + L_2).(H_1 + H_2 - \frac{1}{3} E) = (2 - \frac{1}{9}) - (2 - \frac{1}{3}) = \frac{2}{9}.
\end{align*}
\end{exm} 
\begin{rmk} Naturally, it would be interesting to interpret the right hand side of \eqref{OscillatoryDFIntro} in terms of a suitable stability notion for the mirror Landau-Ginzburg model. Here we only point out a close formal resemblance between \eqref{OscillatoryResidues} and a localised Donaldson-Futaki invariant. Indeed, according to \cite[Section 3]{Legendre_localised}, if $\pi\!:(\mathcal{Y}, [\omega_{\mathcal{Y}}]) \to\PP^1$ is a smooth compactified test configuration for some $n$-dimensional K\"ahler manifold (its general fibre), and if the corresponding $S^1$-action, with infinitesimal generator $v$, has isolated fixed points $Z(v) \cap \mathcal{Y}_0$ on the central fibre, then we have
\begin{align}\label{DFIsolated}
\nonumber\operatorname{DF}(\mathcal{Y}, [\omega_{\mathcal{Y}}]) &= (-1)^{n+1}\sum_{p \in Z(v)\cap \mathcal{Y}_0} \frac{\hat{\mathcal{A}}^n \wedge (\frac{n\mu}{n+1}\hat{\mathcal{A}} -\operatorname{Ric}(\omega_{\mathcal{Y}})-\operatorname{tr} \nabla v + 1)}{\det \nabla v}|_p\\
&=\sum_{p \in Z(v)\cap \mathcal{Y}_0} \frac{\hat{m}^{n} (\frac{n\mu}{n+1}\hat{m} - \operatorname{tr} \nabla v  + 1) }{\det \nabla v} \big|_p, 
\end{align}
where $\hat{\mathcal{A}} = \omega_{\mathcal{Y}} - \hat{m}$, and $\hat{m}$ denotes any Hamiltonian for $v$.

Thus, it seems natural to ask if one can find compactifications $\mathcal{Y}_k$ of $(\C^*)^{n+1}$ and a sequence of K\"ahler forms $\omega_{\mathcal{Y}_k}$ and $S^1$-actions generated by $v_k$, inducing test configurations for compactifications of $(\C^*)^n$, and such that 
\begin{equation*}
\lim_{k\to\infty} \operatorname{DF}(\mathcal{Y}_k, [\omega_{\mathcal{Y}_k}]) = \lim_{k\to\infty}\sum_{p \in \operatorname{Crit}(\mathcal{W}_k)} \frac{ \vartheta^n_{\omega_{\mathcal{X}},k} \big(\frac{n \mu}{n+1}\vartheta_{\omega_{\mathcal{X}},k} - \mathcal{W}_{\operatorname{rel},k}\big)}{\prod_i x^2_i \det \nabla^2 \mathcal{W}_k}\big|_p. 
\end{equation*}
By \eqref{DFIsolated}, when the infinitesimal generators $v_k$ have isolated fixed points on $(\mathcal{Y}_{k})_0$, this means roughly that they should approximate the holomorphic vector fields $\nabla\mathcal{W}_k$, and that the hamiltonians $\hat{m}_k$ should approximate the restriction of the functions $\vartheta_{\omega_{\mathcal{X}}, k}$ to the critical locus $\operatorname{Crit}(\mathcal{W}_k)$.  
\end{rmk}

\subsection{Limiting behaviour}
In the light of \eqref{OscillatoryDFIntro}, it is important to study the behaviour of the residue pairing
\begin{equation*}
\res_{\mathcal{W}(k\omega_{\mathcal{X}})}(\vartheta^n_{\omega_{\mathcal{X}},k}, \frac{n \mu}{n+1}\vartheta_{\omega_{\mathcal{X}},k} - \mathcal{W}_{\operatorname{rel},k})
\end{equation*}
in the large complex structure limit $k \to \infty$. 

In the present work we will study in particular the special case when, roughly speaking, \emph{the mirror of the Donaldson-Futaki invariant concentrates at a single nondegenerate critical point}. We will see that, when this happens, \emph{there is a corresponding simpler constraint on the theta functions on the mirror} (see in particular \eqref{MainIneqIntro}). 

More precisely, we consider the case when there exists a sequence of nondegenerate critical points $p_k \in \operatorname{Crit}(\mathcal{W}(k\omega_{\mathcal{X}}))$ such that  
\begin{align}\label{SpecialLaplaceIntro}
 \operatorname{DF}(\mathcal{X}, [\omega_{\mathcal{X}}]) =  \lim_{k\to\infty} \res_{\mathcal{W}(k\omega_{\mathcal{X}}), p_k}(\vartheta^n_{\omega_{\mathcal{X}},k}, \frac{n \mu}{n+1}\vartheta_{\omega_{\mathcal{X}},k} - \mathcal{W}_{\operatorname{rel},k}),
\end{align}
with the obvious notation for the local residue at a critical point $p_k$. We expect this to happen in some (possibly empty) ``special chamber" of a wall-and-chamber decomposition on the K\"ahler cone of $\mathcal{X}$ (we will actually prove this for toric surfaces in Section \ref{DFAsymp}). Applying \eqref{OscillatoryResidues}, we have in this case
\begin{align*}  
\operatorname{DF}(\mathcal{X}, [\omega_{\mathcal{X}}]) =\lim_{k\to\infty} \frac{ \vartheta^n_{\omega_{\mathcal{X}},k} \big(\frac{n \mu}{n+1}\vartheta_{\omega_{\mathcal{X}},k} - \mathcal{W}_{\operatorname{rel},k}\big)}{\prod_i x^2_i \det \nabla^2 \mathcal{W}_k}\big|_{p_k}. 
\end{align*}
\begin{prop}\label{MainIdentityIntroProp} If $\mathcal{W}(k\omega_{\mathcal{X}})$ has nondegenerate critical points for $k \gg 1$ and the concentration condition \eqref{SpecialLaplaceIntro} holds, then we have
\begin{align}\label{MainIdentityIntro}
 \operatorname{DF}(\mathcal{X}, [\omega_{\mathcal{X}}]) = \lim_{k\to\infty}\vol_{p_k}(\mathcal{X})\left(\frac{n\mu}{n+1} -   \frac{ \mathcal{W}_{\operatorname{rel},k} }{\vartheta_{\omega_{\mathcal{X}},k} }\right)\big|_{p_k}.
\end{align}
In particular, if $\lim_{k\to\infty}\vol_{p_k}(\mathcal{X})$ is real and strictly positive, then $K$-semistability with respect to $(\mathcal{X}, [\omega_{\mathcal{X}}])$ can be expressed as
\begin{equation}\label{MainIneqIntro}
\lim_{k\to\infty} \frac{\mathcal{W}_{\operatorname{rel},k}}{\vartheta_{\omega_{\mathcal{X}},k}}\big|_{p_k} \leq \frac{n\mu}{n+1}.
\end{equation}
\end{prop}
Of course the same result holds if $\lim_{k\to\infty}\vol_{p_k}(\mathcal{X})$ is real and strictly negative, replacing \eqref{MainIneqIntro} with the opposite inequality
\begin{equation*} 
\lim_{k\to\infty} \frac{\mathcal{W}_{\operatorname{rel},k}}{\vartheta_{\omega_{\mathcal{X}},k}}\big|_{p_k} \geq \frac{n\mu}{n+1}.
\end{equation*}   
 
\begin{exm} Our computations in Example \ref{P1ExmCont} show precisely that the concentration condition \eqref{SpecialLaplaceIntro} holds in the basic case of degeneration to the normal cone of a point in $\PP^1$. Our main examples in this paper (covered by Theorem \ref{SlopeThm}) can be seen as higher dimensional generalisations of this fact.

Let us also consider the mirror asymptotics for the volume. Using the formulae of Example \ref{P1Exm}, we compute in this case  
\begin{equation*}
\frac{ \vartheta^{2}_{\omega_{\mathcal{X}},k}}{\prod_i x^2_i \det \nabla^2 \mathcal{W}_k} = \frac{\left(-e^{\frac{8 \pi  k}{3}} x^2 y^2+3 x+3 y\right)^2}{9 x y \left(4-e^{\frac{16 \pi  k}{3}} x^3 y^3\right)},
\end{equation*}
from which, using the asymptotics of critical points worked out in Example \ref{P1ExmCont}, we obtain   
\begin{align*}
& \vol_{p_1}(\cX) \sim -\frac{\left(e^{\frac{2 \pi  k}{3}}+6\right)^2}{9 \left(e^{\frac{4 \pi  k}{3}}-4\right)} \to -\frac{1}{9},\, \vol_{p_i}(\cX) \sim \frac{\left(1\mp6 e^{\frac{\pi  k}{3}}\right)^2}{36 e^{\frac{2 \pi  k}{3}}-9}  \to 1,\,i = 2, 5,\\
& \vol_{p_i}(\cX) \sim -\frac{1}{36 e^{\frac{2 \pi  k}{3}}+9} \to 0,\,i = 3, 4. 
\end{align*}
So we find in this case
\begin{align*}  
\operatorname{DF}(\mathcal{X}, [\omega_{\mathcal{X}}]) &= \lim_{k\to\infty} \vol_{p_1}(\mathcal{X}) \left(1 -  \frac{ \mathcal{W}_{\operatorname{rel},k} }{\vartheta_{\omega_{\mathcal{X}},k} }\right)\big|_{p_1},\,\lim_{k\to\infty} \frac{ \mathcal{W}_{\operatorname{rel},k} }{\vartheta_{\omega_{\mathcal{X}},k} }\big|_{p_1} = 3. 
\end{align*}
\end{exm}
A variant of Proposition \ref{MainIdentityIntroProp} can be obtained by considering a rescaling of the potential 
\begin{equation*}
\cW(\omega_{\cX}) \mapsto c \cW(\omega_{\cX}),\,c\in \C^*.
\end{equation*}
Note that such deformations are contained in the mirror family constructed in \cite{CoatesCortiIritani_hodge}, Section 4.2. The lifts $\mathcal{W}_{\operatorname{rel},k}$, $\vartheta_{\omega_{\mathcal{X}},k}$ can be chosen compatibly with this scaling, that is scaled by the same factor. Each single term in the right hand side of \eqref{OscillatoryResidues} is invariant, while the are nontrivial scalings
\begin{align*}
&\frac{ \vartheta^n_{\omega_{\mathcal{X}},k}}{\prod_i x^2_i \det \nabla^2 \mathcal{W}_k}\big|_p \mapsto \frac{1}{c} \frac{ \vartheta^n_{\omega_{\mathcal{X}},k}}{\prod_i x^2_i \det \nabla^2 \mathcal{W}_k}\big|_p,\\&\big(\frac{n \mu}{n+1}\vartheta_{\omega_{\mathcal{X}},k} - \mathcal{W}_{\operatorname{rel},k}\big)|_p \mapsto c\big(\frac{n \mu}{n+1}\vartheta_{\omega_{\mathcal{X}},k} - \mathcal{W}_{\operatorname{rel},k}\big)|_p.
\end{align*} 
Suppose now that we have
\begin{align*}  
\operatorname{DF}(\mathcal{X}, [\omega_{\mathcal{X}}]) =\lim_{k\to\infty} \frac{ \vartheta^n_{\omega_{\mathcal{X}},k} \big(\frac{n \mu}{n+1}\vartheta_{\omega_{\mathcal{X}},k} - \mathcal{W}_{\operatorname{rel},k}\big)}{\prod_i x^2_i \det \nabla^2 \mathcal{W}_k}\big|_{p_k} < 0.
\end{align*} 
Fix $k$ sufficiently large. Rescaling as above we can assume the normalisation
\begin{equation}\label{ScaledPotential}
\frac{ \vartheta^n_{\omega_{\mathcal{X}},k}}{\prod_i x^2_i \det \nabla^2 \mathcal{W}_k}\big|_{p_k}  = 1.
\end{equation}
\begin{prop}\label{MainIdentityIntroAlternative} Suppose $(\cX, \omega_{\cX})$ is destabilising, and the concentration condition \eqref{SpecialLaplaceIntro} holds. Then, for all sufficiently large $k$, in the scaling for $\cW_k$ given by \eqref{ScaledPotential}, we must have
\begin{equation*}
\frac{n \mu}{n+1}\vartheta_{\omega_{\mathcal{X}},k}(p_k) < \mathcal{W}_{\operatorname{rel},k}(p_k).  
\end{equation*}
\end{prop}
\subsection{Special test configurations}
The geometric meaning of the inequality \eqref{MainIneqIntro} is a little clearer in the Fano case and when $[\omega_{\mathcal{X}}]$ is (semipositive and) a multiple of $c_1(\mathcal{X})$. This is the case for the class of special test configurations in the Fano case, although one needs to allow singularities in general \cite{LiXu}; in the toric Fano case one can take $\mathcal{X}$ smooth by the classical results of \cite{WangZhu_toric}. 
\begin{cor}\label{SpecialTestProp} In the situation of Proposition \ref{MainIdentityIntroProp}, if $[\omega_{\mathcal{X}}]$ is (semipositive and) a multiple of $c_1(\mathcal{X})$, then we have
\begin{equation*}
\operatorname{DF}(\mathcal{X}, [\omega_{\mathcal{X}}]) = \lim_{k\to\infty}\vol_{p_k}(\mathcal{X})\left(\frac{n\mu}{n+1} - \lim_{k\to\infty} \frac{ \mathcal{W}_{\operatorname{rel},k} }{\mathcal{W}_k }\right)\big|_{p_k}.
\end{equation*}
In particular in this case, if $\lim_{k\to\infty}\vol_{p_k}(\mathcal{X})$ is real and strictly positive, the $K$-semistability condition with respect to $(\mathcal{X}, [\omega_{\mathcal{X}}])$ becomes
\begin{align*}
\lim_{k\to \infty}\frac{\mathcal{W}_{k} - \vartheta_{\pi^*\mathcal{O}_{\PP^1}((0) + (\infty)),k} }{\mathcal{W}_{k}}\big|_{p_k}\leq \frac{n}{n+1} \iff \lim_{k\to \infty}\frac{\vartheta_{\mathcal{X}_0,k} +\vartheta_{\mathcal{X}_{\infty},k}}{\mathcal{W}_{k}}\big|_{p_k}\geq \frac{1}{n+1}.  
\end{align*}
Similarly, in the situation of Proposition \ref{MainIdentityIntroAlternative}, if $(\mathcal{X}, [\omega_{\mathcal{X}}])$ is strictly destabilising we must have 
\begin{equation*}
\frac{n \mu}{n+1}\cW_{k}(p_k) < \mathcal{W}_{\operatorname{rel},k}(p_k)
\end{equation*}
for all large $k$.
\end{cor}
In general, we can think of \eqref{MainIneqIntro} as a generalisation of this inequality; for example, when $[\omega_{\mathcal{X}}] = c_1(\mathcal{X}) - \delta [\mathcal{D}_i]$ for some $\mathcal{D}_i$ supported on the central fibre, the inequality becomes
\begin{equation}\label{AnticanInequIntro}
\lim_{k\to \infty} \frac{\mathcal{W}_{k} -\vartheta_{\mathcal{X}_0,k}-\vartheta_{\mathcal{X}_{\infty},k}}{\mathcal{W}_{k} - \delta \vartheta_{\mathcal{D}_i, k}}\big|_{p_k} \leq \frac{n}{n+1}.   
\end{equation}
 
\subsection{Perturbations} The constructions presented so far depend nontrivially on several additional parameters.

On one hand, as we mentioned, the construction of the mirror family given by \cite{CoatesCortiIritani_hodge} is not rigid, but depends on the choice of certain auxiliary parameters, corresponding to unfoldings of the Givental Landau-Ginzburg potential. 

Another important family of perturbations corresponds to allowing angles along the components $D_i \subset D$ in the range $(0, 2\pi)$. Note that in the present paper we only consider angles of $2\pi$.

A more radical possibility would be deforming the stability condition itself, for example, allowing a $B$-field  \cite{ScarpaStoppa_complexified} or generalising $K$-(poly)stability to the stability notions considered in \cite{Dervan_crit_metrics}. 
  
\subsection{Main examples}\label{GeneralExmSec}
In the rest of this paper we work with a surface $X$, either del Pezzo or toric, endowed with a singular, nodal anticanonical divisor $D$. We show how the previous argument can be made fully precise in a number of nontrivial examples, both toric and non-toric. In fact, although we expect that the equivariant theory will be needed in general, our approach in these examples is based on non-equivariant, Givental-type mirror theorems (sometimes partially conjectural), see e.g. \cite{CoatesCortiIritani_hodge, Givental_toric}, and \cite{GrossHackingKeel_LCY} Conjecture 0.19.

Another important simplification for our examples is that we will actually work on the manifold $X$ itself rather than a test configuration $\mathcal{X}$, and, in particular, we will apply the mirror theorems to the lower dimensional manifold $X$.

Such results state that there is an isomorphism of $\C$-algebras
\begin{equation*}
\Theta\!: J(Y, W) \to QH^*(X, [\omega]),
\end{equation*}
the mirror map, from the Jacobi ring 
\begin{equation*}
J(Y, W) = H^0(Y, \mathcal{O}_Y)/(dW) 
\end{equation*}
to the quantum cohomology ring of $(X, [\omega])$, i.e. a suitable ring of formal power series with coefficients in $H^*(X, \mathbb{C})$, endowed with a deformation of the usual cup product, such that 
\begin{equation*}
\Theta([\vartheta_{D_i}]) = [D_i].
\end{equation*}
In particular, we have
\begin{equation*}
\Theta(W) = [-K_X],
\end{equation*}
while the mirror of a K\"ahler class $[\omega] = \sum_i \omega^i [D_i] \in H^{1,1}(X, \R)$ is an element 
\begin{equation*}
[\vartheta_{[\omega]}] = \Theta^{-1}([\omega])= \sum_i \omega^i [\vartheta_{D_i}] \in J(Y, W).
\end{equation*}
Crucially, for our purposes, the mirror map $\Theta$ intertwines the residue pairing on $J(Y, W)$ with the cup product on $H^*(X, \mathbb{C})$. 
 
Suppose now that $(\mathcal{X}, [\omega_{\mathcal{X}}])$ is a K\"ahler test configuration for $(X, [\omega])$, such that the corresponding Donaldson-Futaki invariant $\mathcal{F}$ can be expressed intrinsically on $X$ as an intersection number, 
\begin{equation*}
\mathcal{F}(\mathcal{X}, [\omega_{\mathcal{X}}]) = \sum_{i,j} a_{ij}(\omega)D_i . D_j,\,a_{ij} \in \R. 
\end{equation*}
Examples include suitable product test configurations, degeneration to the normal cone of a divisor $D \subset X$ (leading to the notion of Ross-Thomas slope (in)stability \cite{RossThomas_obstruction}), and more general variants such as flops of degeneration to the normal cone (leading to flop slope (in)stability \cite{CheltsovRubinstein_flops}). Then, by the mirror theorem, we can express the Donaldson-Futaki invariant in terms of theta functions on the mirror as
\begin{equation}\label{mirrorFutakiIntro}
\mathcal{F}(\mathcal{X}, [\omega_{\mathcal{X}}])  = \sum_{i,j} a_{ij} \frac{1}{(2\pi \ii)^2}  \sum_{p\in\operatorname{Crit}(W)}\int_{\Gamma(p)} \frac{\vartheta_{D_i}\vartheta_{D_j}}{\prod_k x_k \del_{x_k} W} \Omega.
\end{equation}
The key problem is how to use such formulae in order to extract an explicit constraint on the theta functions, and so on the complex structure of the mirror $(Y, W)$, from the $K$-(semi)stability constraint $\mathcal{F} \geq 0$. 

As in the general discussion, the main insight we develop here is that $K$-stability is a scale-invariant property, i.e. unchanged by replacing $ \omega  \mapsto  \omega_k:=k\omega$ for $k > 0$, while the mirror map is not. In other words, when $(X, D, [\omega])$ is K-(semi, poly)stable, we get a whole ray of mirror pairs 
\begin{equation*}
(X, D, [\omega_k])\,|\,(Y_k, W_k, \Omega_{Y_k})
\end{equation*}
such that $(X, D, [\omega_k])$ is trivially K-(semi, poly)stable, but $(Y_k, W_k, \Omega_{Y_k})$ changes very nontrivially. In particular, as $k\to \infty$, the variety $Y_k$ degenerates, approaching a so-called large complex structure limit point: the corresponding variety is isomorphic to a union of affine planes (see \cite{GrossHackingKeel_LCY}, Sections 0.5.3, 0.6), 
\begin{equation*}
\mathbb{V}_n \cong \mathbb{A}^2_{12} \cup \cdots \cup \mathbb{A}^2_{n1}.
\end{equation*}
From this viewpoint, we expect that $K$-(semi)stability is a property which can be characterised by a formal neighbourhood of a large complex structure limit of the mirror. 

We study the behaviour of  formulae for Donaldson-Futaki invariants like \eqref{mirrorFutakiIntro} near such large complex structure limits. 

In Sections \ref{AsympSecFano} and \ref{NonToricSec} we show that each critical point gives a well defined contribution in the $k \to \infty$ limit, which is effectively computable (at least in the nondegenerate case). We also spell out the concentration condition \eqref{SpecialLaplaceIntro} explicitly in this case.

We illustrate this method by applying it to product test configurations and to slope (in)stability. 

Suppose $Z = D_i$ is a boundary component. The degeneration to the normal cone $\cX := \Bl_{Z\times\{0\}} X \times \PP^1$, endowed with K\"ahler classes of the form $\pi^*_{X}[\omega] - c [E]$, where $E$ denotes the exceptional divisor, is a test configuration for $(X, [\omega])$ in a natural way. As we recall in Section \ref{SlopeSec}, the Ross-Thomas quotient slope $\mu_{c}(\mathcal{O}_Z, \omega)$ is defined precisely so that $(X, [\omega])$ is ``$K$-destabilised by $Z$", i.e. by $\cX$, precisely when
\begin{equation*} 
\mu_c(\mathcal{O}_Z) < \mu(X).
\end{equation*} 
It is given on the mirror by a complicated expression of the form
\begin{equation*}
\mu_{c}(\mathcal{O}_Z, \omega) = \frac{\sum_{p \in \operatorname{Crit}(W_k)} d_p(\omega_k)}{\sum_{p \in \operatorname{Crit}(W_k)} r_p(\omega_k)},
\end{equation*}
valid for all $k >0$ (see Section \ref{SlopeSec} for the details). However we show that, in many (nondegenerate) cases, there is a simpler expression nearby a large complex structure limit, namely there exists a sequence of critical points $p_k$ for $W_k$ such that
\begin{equation}\label{slopeSymmetry}
\mu_{c}(\mathcal{O}_Z, \omega) = \lim_{k\to \infty} \frac{d_{p_k}(\omega_k)}{r_{p_k}(\omega_k)}.
\end{equation}
This corresponds to the concentration condition \eqref{SpecialLaplaceIntro}, and we give several explicit examples where this holds. In all these cases, the following result (which is indeed a version of the identities \eqref{SpecialLaplaceIntro}, \eqref{MainIdentityIntro} in this special case) can be applied effectively.  
\begin{thm}\label{SlopeThm} For any sequence of critical points $p_k$ satisfying \eqref{slopeSymmetry}, the quotient slope can be computed in terms of theta functions on the mirror as  
\begin{equation*}
\mu_{c}(\mathcal{O}_Z) = \frac{1}{c}\lim_{k\to \infty} \frac{ 1 - 3 c\left(\frac{k \vartheta_Z }{\vartheta_{[\omega_k]} }-\frac{k W_k }{\vartheta_{[\omega_k]} }\right)}{ 1  -2 c \frac{k\vartheta_Z }{\vartheta_{[\omega_k]} } }\big|_{p_k}.
\end{equation*}
This holds for suitable K\"ahler classes on toric del Pezzos, their simple degenerations, and the non-toric del Pezzo surfaces of degree $5$ and $4$. In particular, in the case of the anticanonical polarisation $[\omega] = -K_X$, we have
\begin{equation*}
\mu_{c}(\mathcal{O}_Z) = \frac{1}{c}\lim_{k\to \infty} \frac{ 1 - 3 c\left(\frac{\vartheta_Z }{W_k }-1\right)}{ 1  -2 c \frac{\vartheta_Z }{W_k } }\big|_{p_k}.
\end{equation*}
\end{thm}
\begin{rmk}\label{ThimbleRmk} Suppose $Z$ destabilises $(X, [\omega])$. Then, for some fixed $k \gg 1$, we have
\begin{equation*}
\frac{1}{c} \operatorname{Re}\frac{ 1 - 3 c\left(\frac{k \vartheta_Z }{\vartheta_{[\omega_k]} }-\frac{k W_k }{\vartheta_{[\omega_k]} }\right)}{ 1  -2 c \frac{k\vartheta_Z }{\vartheta_{[\omega_k]} } }\big|_{p_k} < \mu(X).
\end{equation*}
The critical point $p_k$ corresponds to a Lagrangian thimble $L_k(Z)$ for $W_k \!: Y_k \to \C$ (see \cite{Auroux_delPezzo}, Section 4.1). This suggests that $(Y_k, W_k, \Omega_{Y_k})$ is ``destabilised" by the Lagrangian submanifold $L_k(Z)$.
\end{rmk}
Theorem \ref{SlopeThm} is proved in Section \ref{SlopeSec} using our general results in Sections \ref{AsympSecFano} and \ref{NonToricSec}: in particular our discussion of  toric del Pezzos (see Sections \ref{BlpSlopeSubsec}-\ref{BlpqSlopeSubsec}), their simple degenerations (Section \ref{IteratedBlpSlopeSec}), and the non-toric del Pezzo surfaces of degree $5$ (\ref{Deg5AsympSec}-\ref{Deg5SlopeSec}) and degree $4$ (\ref{Deg4AsympSec}-\ref{Deg4SlopeSec}).

Thus we also have a precise version of the inequality \eqref{AnticanInequIntro},  
in cases when $[\omega]$ is a multiple of $c_1(X)$. 
\begin{cor} $(X, -K_X)$ is slope semistable with respect to a divisor $Z \subset X$ as in Theorem \ref{SlopeThm}, iff on the mirror we have 
\begin{equation*}
\lim_{k\to\infty}\frac{(3 + c^{-1})W_k-3\vartheta_Z}{W_k - 2 c \vartheta_Z}\big|_{p_k} \geq 1.
\end{equation*}
\end{cor}
This can be used to check (in)stability with respect to a boundary component $Z$ in some of our examples, see Sections \ref{BlpqSlopeSubsec} and \ref{IteratedBlpSlopeSec}.
 
\subsection{Product configurations} We also discuss a variant of these results in the case of product test configurations, i.e. Hamiltonian holomorphic vector fields $V_a$, with normalised Hamiltonian $\mu_a$. As an example, in Section \ref{ProductTestSec} we show the following result and provide an example where it can be applied effectively.
\begin{thm} Suppose that the Futaki invariant $\mathcal{F}_{[\omega]}(J V_a)$ vanishes, and the only contributions to its localisation come from two $V_a$-fixed divisors $D_1, D_2$ (in particular, the contributions from isolated points cancel out). Then, if $\mathcal{F}_{[\omega]}(J V_a)$ concentrates uniformly on sets of critical points $\{p_k\}, \{p'_k\}$ in the large complex structure limit, we must have
 \begin{align*}
& \lim_{k\to \infty} \kappa(p_k,p'_k)  \frac{\vartheta_{D_1}\big(1-\pi \frac{\mu_a(D_1)}{\langle \operatorname{w}^{D_1}, a\rangle} (\frac{W_k}{\vartheta_{[\omega_k]}} + \frac{\vartheta_{D_1}}{\vartheta_{[\omega_k]}})\big)|_{p_k}}{\vartheta_{D_2}\big(1-\pi \frac{\mu_a(D_2)}{\langle \operatorname{w}^{D_2}, a\rangle} (\frac{W_k}{\vartheta_{[\omega_k]} } + \frac{\vartheta_{D_2} }{\vartheta_{[\omega_k]} })\big)|_{p'_k}} = \frac{\#\{p'\}}{\#\{p\}} \frac{\mu_a(D_2)}{\mu_a(D_1)},
 \end{align*} 
where for any pair of critical points $p, p'$ we set
\begin{equation*}
\kappa(p,p') = \frac{\vartheta_{[\omega_k]}(p')}{\vartheta_{[\omega_k]}(p)} \frac{\vol_{p,k}(\omega)}{\vol_{p', k}(\omega)}.
\end{equation*}
In particular, for the anticanonical polarisation we have
 \begin{align*}
& \lim_{k\to \infty} \kappa(p_k,p'_k) \frac{\vartheta_{D_1}\big(1-\pi(1 + \frac{\vartheta_{D_1}(p)}{W_k(p)})\big)|_{p_k}}{\vartheta_{D_2}\big(1-\pi(1 + \frac{\vartheta_{D_2}(p')}{W_k(p')})\big)|_{p'_k}} = \frac{\#\{p'\}}{\#\{p\}} \frac{\mu_a(D_2)}{\mu_a(D_1)} = \frac{\#\{p'\}}{\#\{p\}} \frac{\langle \operatorname{w}^{D_2}, a\rangle}{\langle \operatorname{w}^{D_1}, a\rangle},
 \end{align*} 
 where
 \begin{equation*}
\kappa(p,p') = \frac{W_k(p')}{W_k(p)} \frac{\vol_{p,k}(-K_X)}{\vol_{p', k}(-K_X)}.
\end{equation*}
\end{thm}

\section{Toric del Pezzo case}\label{AsympSecFano}
Suppose $(X, D = D_1 + \dots +D_n)$ is given by a toric del Pezzo with its toric boundary. Write $\{v_i,\,i=1,\dots,n\}$ for the generators of the fan of $X$, with $v_i$ dual to the divisor $D_i$. Then, as explained in Example \ref{GiventalFormula}, we have explicit formulae for theta functions (see e.g. \cite{CoatesCortiIritani_hodge, Givental_toric})
\begin{align*}
&\vartheta_{D_i}(\omega) = a_i(\omega) x^{v_i},\,W(\omega) = \sum^{n}_{i=1} \vartheta_{D_i}(\omega).
\end{align*}
Mirror symmetry shows in particular that we have an identity between intersection numbers and residue pairings
\begin{equation*}
D_i . D_j = \frac{1}{(2\pi \ii)^2} \sum_{p\in\operatorname{Crit}(W)}\int_{\Gamma(p)} \frac{\vartheta_{D_i}\vartheta_{D_j}}{\prod_k x_k \del_{x_k} W} \Omega,
\end{equation*} 
where $\Gamma(p) \subset Y$ is a compact integration cycle around a critical point $p$ of $W$, and $\Omega$ denotes the holomorphic volume form, given by 
\begin{equation*}
\Omega = \prod_k \frac{d x_k}{x_k}
\end{equation*}
in the toric del Pezzo case. According to \cite{CoatesCortiIritani_hodge}, Section 6, we can also express the contribution of a \emph{nondegenerate} critical point to the residue pairing as
\begin{equation*}
\int_{\Gamma(p)} \frac{\vartheta_{D_i}\vartheta_{D_j}}{\prod_k x_k \del_{x_k} W} \Omega = \frac{\vartheta_{D_i}(p)\vartheta_{D_j}(p)}{(p_1 p_2)^2 \det \nabla^2 W(p)},
\end{equation*}
where we write $p = (p_1, p_2)$ for the torus coordinates of $p$ in our case. 

\subsection{Asymptotics of critical points}
We now rescale the K\"ahler form by $\omega_k = k \omega$, $k > 0$, and study the $k \to \infty$ limit of
\begin{equation*}
\vartheta_{D_i, k} = \vartheta_{D_i}(\omega_k),\, W_k = W(\omega_k).  
\end{equation*}
The general theory developed in \cite{CoatesCortiIritani_hodge}, Section 6.1, shows that there exist expansions
\begin{equation*}
p_{j} = \alpha_{p, j} e^{2\pi k \beta_{p,j}}(1+O(k^{-1})) \text{ for } p\in\operatorname{Crit}(W_k)
\end{equation*}
for $\alpha_{p, j} \in \C^*,\,\beta_{p,j} \in \R,\,j=1,2,\,k\to \infty$. Write $(x, y) = (x_1, x_2)$ for our torus coordinates of a general point. Note that, by the defining property, $a_i(k\omega) = (a_i(\omega))^k$. Then, we have
\begin{equation*}
x_j \del_{x_j} W_k = \sum^n_{i=1} v_{i, j} e^{k \log a_i(\omega)}  x^{v_i},\,j=1,2, 
\end{equation*} 
and so, at a critical point $p$ of $W_k$,  
\begin{equation*}
\sum^n_{i=1} v_{i, j} e^{k\log a_i(\omega)}  p^{v_i,1}_1 p^{v_i,2}_2 = 0,\,j=1,2.
\end{equation*} 
By the existence of the expansion, we must have, for $p\in\operatorname{Crit}(W_k)$,
\begin{align*}
\sum_{i\,:\,v_{i, j} \neq 0,\,\langle v_i , \beta_p \rangle- \frac{1}{2\pi}\log a_i(\omega) \text{ is maximal}} v_{i, j} (\alpha_{p, 1})^{v_{i,1}} (\alpha_{p, 2})^{v_{i,2}}  = 0,\,j=1,2.
\end{align*} 
This is a constraint which allows to determine both $\beta_p$ and $\alpha_p$, given the existence results. Namely, the set $\{i\,:\,v_{i,j}\neq 0,\,\langle v_i , \beta_p \rangle- \frac{1}{2\pi}\log a_i(\omega) \text{ is maximal}\}$ must contain at least two distinct elements $i_j, i'_j$: if this consisted of a single $i_j$, the corresponding equation for the coefficients would be
\begin{equation*}
v_{i_j,j} (\alpha_{p, 1})^{v_{i_j, 1}} (\alpha_{p, 2})^{v_{i_j,2}}  = 0,
\end{equation*}  
which does not have solutions on the torus $(\C^*)^2$. Thus, we obtain the system of affine linear inequalities for $\beta_p$
\begin{align*}
&\langle v_{i_j} , \beta_p\rangle - \frac{1}{2\pi}\log a_{i_j}(\omega)= \langle v_{i'_j} , \beta_p \rangle - \frac{1}{2\pi}\log a_{i'_j}(\omega)\\&\geq \langle v_{i} , \beta_p \rangle - \frac{1}{2\pi}\log a_i(\omega),\,j=1,2, i \neq i_{j}, i'_j. 
\end{align*}  
This defines a wall-and-chamber structure on the space of K\"ahler parameters $[\omega]$, corresponding to the possible asymptotics of the critical points. Then, the coefficients $\alpha_p$ are determined by the algebraic equations above. 
\begin{exm} Consider $X = \operatorname{Bl}_p \PP^2$ with K\"ahler class $\omega_k = k (H- q E)$, $q \in (0,1)$. A Landau-Ginzburg potential is given by
\begin{align*}
&W_k =  \frac{e^{-2\pi k }}{x y} +   x +  y + e^{2\pi k q } x y 
\end{align*}
(this can be modified by scaling the variables $x$, $y$). Let us determine the wall-and-chamber structure on the K\"ahler parameters determined by the asymptotics of critical points, i.e. the corresponding decomposition of the interval $(0,1)$. The critical locus is cut out by the equations
\begin{align*}
&\del_x W_k = y e^{2 \pi  k q}-\frac{e^{-2 \pi  k}}{x^2 y}+1 = 0,\,\del_y W_k = x e^{2 \pi  k q}-\frac{e^{-2 \pi  k}}{x y^2}+1 = 0.
\end{align*}
Using the asymptotics $x = a e^{2\pi k b}(1 + O(k^{-1}))$, $y = c e^{2\pi k d}(1 + O(k^{-1}))$ gives the conditions
\begin{align*}
& -\frac{e^{-4 \pi  b k-2 \pi  d k-2 \pi  k}}{a^2 c}+c e^{2 \pi  d k+2 \pi  k q}+1 = 0,\,-\frac{e^{-2 \pi  b k-4 \pi  d k-2 \pi  k}}{a c^2}+a e^{2 \pi  b k+2 \pi  k q}+1 = 0, 
\end{align*}
up to terms vanishing as $k\to \infty$. Following the general procedure we described, consider the conditions
\begin{align*}
2 \pi  d k+2 \pi  k q = 0,\,2 \pi  b k+2 \pi  k q = 0 \iff b = d = -q. 
\end{align*}
\end{exm}
The critical equations become
\begin{align*}
-\frac{e^{2 \pi  k (3 q-1)}}{a^2 c}+c+1 = 0,\,-\frac{e^{2 \pi  k (3 q-1)}}{a c^2}+a+1=0,
\end{align*}
so, for $q \in (0, \frac{1}{3})$, there is a critical point with asymptotics $(x, y) \sim (-e^{2\pi k q}, -e^{2\pi k q})$. Similarly, the conditions
\begin{align*}
 -4 \pi  b k-2 \pi  d k-2 \pi  k = 0,\,-2 \pi  b k-4 \pi  d k-2 \pi  k =0 \iff b = d = -\frac{1}{3}
\end{align*} 
give critical equations
\begin{equation*}
-\frac{1}{a^2 c}+c e^{\frac{2}{3} \pi  k (3 q-1)}+1 = 0,\,-\frac{1}{a c^2}+a e^{\frac{2}{3} \pi  k (3 q-1)}+1 =0, 
\end{equation*}
so, for $q \in (0, \frac{1}{3})$, the remaining three critical points satisfy $(x, y) \sim (\xi e^{-\frac{2}{3}\pi k}, \xi e^{-\frac{2}{3}\pi k})$ for all roots $\xi$ of $a^3 - 1 = 0$.  

There is a wall at $q = \frac{1}{3}$ (corresponding to the anticanonical polarisation), for which the critical points satisfy $(x, y) \sim (\xi e^{-\frac{2}{3}\pi k}, \xi e^{-\frac{2}{3}\pi k})$ for all roots $\xi$ of $a^4 + a^3 - 1 = 0$. 

Finally, for $q \in (\frac{1}{3}, 1)$, one checks that the critical asymptotics are determined by the conditions
\begin{equation*}
-4 \pi  b k-2 \pi  d k-2 \pi  k = 2 \pi  d k+2 \pi  k q = 2 \pi  b k+2 \pi  k q \iff b = d = \frac{1}{4}(-1-q)
\end{equation*} 
giving critical equations
\begin{equation*}
-\frac{e^{\frac{1}{2} \pi  k (3 q-1)}}{a^3}+a e^{\frac{1}{2} \pi  k (3 q-1)}+1 =0,\,-\frac{e^{\frac{1}{2} \pi  k (3 q-1)}}{a^3}+a e^{\frac{1}{2} \pi  k (3 q-1)}+1=0.
\end{equation*}
It follows that we must have $(x, y) \sim (\xi e^{-\frac{\pi k}{2}(1+q)}, \xi e^{-\frac{\pi k}{2}(1+q)})$ where $\xi$ is any root of $a^4 - 1 =0$. 
\subsection{Asymptotics of residues}
Differentiating the expression
\begin{equation*}
x_r \del_{x_r} W_k = \sum^n_{i=1} v_{i, r} e^{k \log a_i(\omega)}  x^{v_i}  
\end{equation*}
gives
\begin{equation*}
x_s \del_{x_s}\big(x_r \del_{x_r} W_k\big) = \sum^n_{i=1} v_{i, r} v_{i, s} e^{k \log a_i(\omega)}  x^{v_i},  
\end{equation*}
and so
\begin{align*}
& x_r x_s \del^2_{x_r x_s} W_k = \sum^n_{i=1} v_{i, r} v_{i, s} e^{k \log a_i(\omega)}  x^{v_i} - \delta_{rs} x_s \del_{x_s} W_k\\
& = \sum^n_{i=1} (v_{i, r} v_{i, s} -\delta_{rs} v_{i,s}) e^{k \log a_i(\omega)}  x^{v_i} 
\end{align*}
($\delta_{rs}$ denoting the Kronecker delta). So we have
\begin{align*}
&(x y)^2 \det \nabla^2 W_k = (x^2\del^2_{x,x} W_k )(y^2\del^2_{y,y} W_k) - (x y \del_{x,y} W_k)^2\\
&= \sum^n_{i, j = 1} (v^2_{i,1} v^2_{j,2} - v^2_{i,1} v_{j,2} - v_{i, 1} v^2_{j,2}) e^{k \log a_i(\omega)a_j(\omega)} x^{v_i + v_j}. 
\end{align*}
Using the asymptotic expansion for critical points  
\begin{equation*}
p_{r} = \alpha_{p, r} e^{2\pi k \beta_{p,r}}(1+O(k^{-1})) \text{ for } p\in\operatorname{Crit}(W_k),\,\alpha_{p, r} \neq 0,\,r=1,2,\,k\to \infty 
\end{equation*}
we find
\begin{align*}
&(x y)^2 \det \nabla^2 W_k\big|_p \\
&= \sum^n_{i, j = 1} (v^2_{i,1} v^2_{j,2} - v^2_{i,1} v_{j,2} - v_{i, 1} v^2_{j,2}) e^{k \log a_i(\omega)a_j(\omega)} \alpha^{v_{i,1} + v_{j,1}}_{p, 1} \alpha^{v_{i,2} + v_{j,2}}_{p, 2} e^{2\pi k \langle (v_i + v_j), \beta_p\rangle} (1 + O(k^{-1})), 
\end{align*}
and so we have, for a nondegenerate critical point,
\begin{align*}
&(x y)^{-2} (\det \nabla^2 W_k)^{-1}\big|_p\\
&=\big( \sum_{i, j \,:\, \langle  v_i + v_j , \beta_p\rangle -\frac{1}{2\pi} \log a_ia_j\text{ is maximal}} e^{-k \log a_i(\omega)a_j(\omega)}e^{-2\pi k \langle (v_i + v_j), \beta_p\rangle}\\
&(v^2_{i,1} v^2_{j,2} - v^2_{i,1} v_{j,2} - v_{i, 1} v^2_{j,2})  \alpha^{v_{i,1} + v_{j,1}}_{p, 1} \alpha^{v_{i,2} + v_{j,2}}_{p, 2}  \big)^{-1}(1 + O(k^{-1})).
\end{align*}
Similarly, for theta functions, we have
\begin{equation*}
\vartheta_{D_r, k} := \vartheta_{D_r}(k \omega) = e^{k \log a_r(\omega)} x^{v_i},
\end{equation*}
and so
\begin{equation*}
\vartheta_{D_r, k}|_p = (\alpha_{r,1})^{v_{r,1}} (\alpha_{r, 2})^{v_{r, 2}} e^{k \log a_r(\omega)}e^{2\pi k \langle v_r,\beta_p\rangle}(1+O(k^{-1})).
\end{equation*}
It follows that the contribution of a critical point $p$ to the residue pairing satisfies
\begin{align}\label{ResidueAsymp}
\nonumber \frac{\vartheta_{D_r,k} \vartheta_{D_s,k} }{(x y)^2 \det \nabla^2 W }\big|_p
&= \big( \sum_{i, j \,:\, \langle  v_i + v_j , \beta_p\rangle -\frac{1}{2\pi} \log a_ia_j\text{ is maximal}} e^{-k \log a_i(\omega)a_j(\omega)}e^{-2\pi k \langle (v_i + v_j), \beta_p\rangle}\\
\nonumber&(v^2_{i,1} v^2_{j,2} - v^2_{i,1} v_{j,2} - v_{i, 1} v^2_{j,2})  \alpha^{v_{i,1} + v_{j,1}}_{p, 1} \alpha^{v_{i,2} + v_{j,2}}_{p, 2}  \big)^{-1}\\
&\cdot e^{k \log a_r(\omega)a_s(\omega)} \alpha^{v_r + v_s}_p e^{2\pi k \langle v_r + v_s,\beta_p\rangle}(1+O(k^{-1})).
\end{align}
\subsection{Donaldson-Futaki invariants}\label{DFAsymp}
Let $(\mathcal{X}, [\omega_{\mathcal{X}}])$ be a K\"ahler test configuration for $(X, [\omega])$, such that the corresponding Donaldson-Futaki invariants can be expressed \emph{as an intersection number on $X$}, 
\begin{equation*}
\mathcal{F}(\mathcal{X}, [\omega_{\mathcal{X}}]) = \sum_{r, s} a_{rs}(\omega_{\mathcal{X}}) D_r.D_s.
\end{equation*}
Applying the mirror theorem, if all critical points are nondegenerate, we obtain
\begin{equation*}
\mathcal{F}(\mathcal{X}, [\omega_{\mathcal{X}}]) =  \sum_{r,s} a_{rs}( \omega_{\mathcal{X}})\sum_{p \in \operatorname{Crit}(W_k)}\frac{\vartheta_{D_r,k} \vartheta_{D_s,k} }{(x y)^2 \det \nabla^2 W_k }\big|_p.
\end{equation*}
Note that we have
\begin{equation*}
\frac{\vartheta_{D_r,k} \vartheta_{D_s,k} }{(x y)^2 \det \nabla^2 W_k }|_{p_k} = c_{r s} e^{2\pi k\gamma_{r s}}(1+O(k^{-1}))
\end{equation*}
for unique $c_{rs}, \gamma_{rs}$ determined by \eqref{ResidueAsymp}, which yields an expansion  
\begin{align*}
&\sum_{r,s}( a_{r s}( \omega_{\mathcal{X}}))\frac{\vartheta_{D_r,k} \vartheta_{D_s,k} }{(x y)^2 \det \nabla^2 W_k }|_{p_k} =  c^a(p_k) e^{2\pi k\gamma^a(p_k)}(1+O(k^{-1})) 
\end{align*} 
for unique $c^a$, $\gamma^a$, computed by \eqref{ResidueAsymp}. So, we have
\begin{equation*}
\mathcal{F}(\mathcal{X}, [\omega_{\mathcal{X}}]) =  \lim_{k\to \infty}  \sum_{p_k\in \operatorname{Crit}(W_k)\,:\,\gamma^a(p_k) = 0} c^a(p_k).
\end{equation*}
In particular, in the special case when 
\begin{equation*} 
\#\{p_k\!: \gamma^a(p_k) = 0\} = 1 
\end{equation*}
we have
\begin{equation*}
 \mathcal{F}(\mathcal{X}, [\omega_{\mathcal{X}}]) = \lim_{k\to \infty} c^a(\hat{p}_k),
\end{equation*}
where $\{p_k\!: \gamma^a(p_k) = 0\} = \{\hat{p}_k\}$. This is a version of the identity \eqref{SpecialLaplaceIntro}, implying \eqref{MainIdentityIntro}, in our situation, and indeed by \eqref{ResidueAsymp} we see that it holds in some (possibly empty) chamber in the space of K\"ahler classes $[\omega_{\mathcal{X}}]$ (in particular, possibly varying the K\"ahler class $[\omega]$).
\subsection{Toric weak del Pezzos} 
When $(X, D)$ is a toric weak del Pezzo surface with its toric boundary the mirror map is much more complicated (it is constructed explicitly as a special case in \cite{CoatesCortiIritani_hodge}, Section 4.6; see also \cite{Chan_toric}), but we still have
\begin{equation*}
W(\omega) = \sum^n_{i = 1} \vartheta_{D_i}(\omega), 
\end{equation*}
where now
\begin{equation*}
\vartheta_{D_i}(\omega) = e^{k \log a_i(\omega)}(1 + O(k^{-1})) x^{v_i},\,k\to\infty.
\end{equation*}
So, the analysis carried out in the toric del Pezzo case remains valid up to $O(k^{-1})$ terms in the large complex structure limit $k \to \infty$. 

\section{Non-toric case}\label{NonToricSec}
Mirror Landau-Ginzburg models for del Pezzo surfaces endowed with a maximally degenerate anticanonical divisor have been constructed in \cite{GrossHackingKeel_LCY}. The mirror theorem in this case is partially conjectural. Explicit equations for the mirror Landau-Ginzburg models have been studied further e.g. in \cite{Arguz_equations, Barrott_equations}. Here we restrict to two examples, the del Pezzo surfaces of degrees $5$ and $4$, for which we compute the asymptotic behaviour of the critical points of the LG potentials in the large complex structure limit, with respect to certain K\"ahler parameters. This is later used to complete the proof of Theorem \ref{SlopeThm} in these cases. 
\subsection{Degree $5$ del Pezzo}\label{Deg5AsympSec}
Let $(X, D = D_1 + \cdots D_5)$ denote a degree $5$ del Pezzo,
\begin{equation*}
X \cong \operatorname{Bl}_{\{p_i\}}\PP^2,\,i=1,\dots4,\,\{p_i\} \text{ generic}, 
\end{equation*}
together with a fixed cycle of smooth rational $-1$ curves $D$. Writing $E_i$, $L_{ij}$ for the obvious exceptional divisors and proper transforms of lines, we can choose 
\begin{equation*}
D_1 = L_{12},\,D_2 = E_2,\,D_3 = L_{23},\, D_4 = E_3,\,D_5 = L_{34}.
\end{equation*} 
Note that, for each $i = 1,\dots, 4$, there is a unique $-1$ curve $\tilde{E}_i$ which is not contained in $D$ and intersects $D_i$ transversely, namely
\begin{equation*}
\tilde{E}_1 = E_1,\,\tilde{E}_2 = L_{24},\,\tilde{E}_3 = L_{14},\,\tilde{E_4} = L_{13} ,\,\tilde{E}_5 = E_4.
\end{equation*}
According to \cite{GrossHackingKeel_LCY}, the mirror family has underlying affine variety 
\begin{equation*}
\mathcal{Y} \subset \operatorname{Spec}\C[\vartheta_{1}, \cdots,\vartheta_5] \times \operatorname{Spec} \C[\operatorname{NE}(X)]
\end{equation*}
cut out by the equations
\begin{equation*}
\vartheta_{i-1} \vartheta_{i+1} = z^{[D_i]}\big(\vartheta_{i} + z^{[\tilde{E}_i]}\big),\,i=1,\dots, 5 
\end{equation*}
(cyclic ordering). So, the mirror to $(X, D)$ endowed with a K\"ahler class $[\omega]$ is given by the affine surface
\begin{equation*}
Y_{\omega} = \{\vartheta_{i-1} \vartheta_{i+1} = e^{-2\pi \int_{D_i}\omega}\big(\vartheta_{i} + e^{-2\pi\int_{\tilde{E}_i}\omega}\big)\} \subset \operatorname{Spec}\C[\vartheta_{1}, \cdots,\vartheta_5],
\end{equation*}
together with the Landau-Ginzburg potential
\begin{equation*}
W_{\omega} = \big(\sum^5_{i=1}\vartheta_i\big) |_{Y_{\omega}}\!: Y_{\omega} \to \C.
\end{equation*}
According to \cite{GrossHackingKeel_LCY}, Example 3.7, there is a maximal dense open subset $\mathcal{U} \subset Y$ where $\vartheta_1, \vartheta_{2}$ give holomorphic local coordinates $x := \vartheta_1, y := \vartheta_2$, for which
\begin{equation*}
\Omega = \frac{d x}{x} \wedge \frac{d y}{y}.
\end{equation*}
We compute, on $\mathcal{U}$,
\begin{align*}
&\vartheta_3 = z^{[D_2 + \tilde{E}_2]}\vartheta^{-1}_1 + z^{[D_2]}\vartheta^{-1}_1\vartheta_2 = z^{[H - E_4]}\vartheta^{-1}_1 + z^{[E_2]}\vartheta^{-1}_1\vartheta_2,\\
&\vartheta_4 = z^{[D_3]} \vartheta^{-1}_2 \vartheta_3 + z^{[D_3 + \tilde{E}_3]}\vartheta^{-1}_2 \\
&= z^{[H-E_3]}\vartheta^{-1}_1 + z^{[2H-\sum^4_{i=1}E_i]}\vartheta^{-1}_2 + z^{[2H-E_2-E_3-E_4]}\vartheta^{-1}_1\vartheta^{-1}_2,\\
&\vartheta_5 = z^{[D_1 + \tilde{E}_1]}\vartheta^{-1}_2 + z^{[D_1]}\vartheta_1\vartheta^{-1}_2  = z^{[H-E_2]}\vartheta^{-1}_2 + z^{[H-E_1-E_2]}\vartheta_1\vartheta^{-1}_2.
 \end{align*}
So, fixing a K\"ahler class 
\begin{equation*}
[\omega] = H - \sum^4_{i=1} a_i  E_i,
\end{equation*}
we have
\begin{align*}
&\vartheta_3 = e^{-2\pi(1-a_4)}\vartheta^{-1}_1 + e^{-2\pi a_2 }\vartheta^{-1}_1\vartheta_2,\\
&\vartheta_4  = e^{-2\pi(1-a_3)}\vartheta^{-1}_1 + e^{-2\pi (2-\sum^4_{i=1}a_i)}\vartheta^{-1}_2 + e^{-2\pi(2 -a_2-a_3-a_4)}\vartheta^{-1}_1\vartheta^{-1}_2\\
&\vartheta_5 = e^{-2\pi(1-a_2)}\vartheta^{-1}_2 + e^{-2\pi(1- a_1- a_2)}\vartheta_1\vartheta^{-1}_2.
 \end{align*}
Let us consider the case of 
\begin{equation*}
[\omega] = [H - \frac{1}{2} E_1 - \delta \sum^4_{i=2} E_i],
\end{equation*}
with corresponding potential $W^{\delta}$. Then, for $\delta > 0$ sufficiently small, the asymptotic behaviour as $k\to \infty$ of the points of $\operatorname{Crit}(W^{\delta})\cap \mathcal{U}$ converges to that of the critical points of the limiting potential $W^0|_{\mathcal{U}}\!:(\C^*)^2 \to \C$, given by
\begin{equation*}
W^0|_{\mathcal{U}} = \frac{e^{-k\pi} x}{y}+\frac{e^{-4 \pi  k}}{x y}+\frac{2 e^{-2 \pi 
   k}}{x}+\frac{e^{-3 \pi  k}}{y}+\frac{e^{-2 \pi  k}}{y}+\frac{y}{x}+x+y,
\end{equation*}
By direct computation, we have
\begin{align*}
& \del_x W^0 = 0 \iff -e^{-\pi k} x^2+2 e^{-2 \pi  k} y+e^{-4 \pi  k}-x^2 y+y^2 = 0, \\
& \del_y W^0 = 0 \iff e^{-\pi k} x^2+e^{-3 \pi  k} x+e^{-2 \pi  k} x+e^{-4 \pi  k}-(x+1) y^2 = 0,
\end{align*}
and solving for $x$, respectively $y$ shows that at points $(x, y) \in \operatorname{Crit}(W_{0})\cap\mathcal{U}$ we have the relations
\begin{align*}
 x = \pm \frac{e^{2 \pi  k} y+1}{\sqrt{e^{4 \pi  k} y+e^{3 \pi  k}}},\, y = \pm \frac{\sqrt{e^{3 \pi  k} x^2+e^{\pi  k} x+e^{2 \pi  k} x+1}}{\sqrt{e^{4 \pi  k} (x+1)}}.
\end{align*}
We claim that these relations imply that that the possible asymptotics are 
\begin{equation*}
x \sim -1,\,y\sim1,\, x\sim a_{i,1} e^{2\pi k \beta_{i,1}}, y \sim a_{i,2}e^{2\pi k \beta_{i, 2}},\, \beta_{i,j } < 0.
\end{equation*}
Indeed, if we have $y \sim a_{i,2} e^{2\pi k \beta_{i, 2}}$ for $\beta_{i, 2} > 0$, then the relation for $x$ implies $x \sim a_{i,1} e^{\pi k \beta_{i,2}}$, and the relation for $y$ gives in turn the constraint
\begin{equation*}
a_{i,2} e^{2\pi k \beta_{i,2}} \sim a_{i,2} e^{\pi k (\frac{3}{2} +2\beta_{i,2}-2-\frac{\beta_{i,2}}{2})} \iff \beta_{i,2} = -1, 
\end{equation*}
a contradiction. On the other hand, the compatibility of $x \sim -1,\,y\sim1$ is easily checked. 
\subsection{Degree $4$ del Pezzo}\label{Deg4AsympSec} 
Following \cite{Arguz_equations, Barrott_equations}, we consider the case when $X = \operatorname{Bl}_{\{p_i\}}\PP^2$, $i =1,\cdots, 5$, $\{p_i\}$ generic, endowed with the anticanonical cycle of $-1$ curves $D = D_1 + D_2 + D_3 + D_4$, such that (in the usual notation)
\begin{equation*}
D_1 = E_1,\,D_2 = L_{12},\,D_3 = L_{34},\,D_4 = L_{15}. 
\end{equation*}
Then, the mirror Landau-Ginzburg family has total space
\begin{equation*}
\mathcal{Y} \subset \operatorname{Spec}\C[\vartheta_{1}, \cdots,\vartheta_4] \times \operatorname{Spec} \C[\operatorname{NE}(X)]
\end{equation*}
cut out by the equations
\begin{align*}
&\vartheta_{1}\vartheta_{3} = z^{[D_2]}\vartheta_2 + z^{[D_4]}\vartheta_4+\\
&z^{[H-E_1]}+z^{[2H-E_1-E_2-E_3-E_5]} + z^{[2H-E_1-E_2-E_4-E_5]},\\
&\vartheta_{2}\vartheta_{4} = z^{[D_1]}\vartheta_1 + z^{[D_3]}\vartheta_3+\\
&z^{[H-E_3]}+z^{[ H-E_4]} + z^{[2H-E_2-E_3-E_4-E_5]},
\end{align*}
and Landau-Ginzburg potential $W = \sum^4_{i=1}\vartheta_i$. According to \cite{Arguz_equations, Barrott_equations}, there is a dense open subset $\mathcal{U}\subset \mathcal{Y}$, biholomorphic to $(\C^*)^2$ with torus coordinates $x, y$, such that on $\mathcal{U}$ we have
\begin{align*}
&\vartheta_{D_1} = \frac{z^{[E_1-E_5]}}{x y^2}+\frac{z^{[H-E_4-E_5]}}{y}+\frac{1}{x y},\\
&\vartheta_{D_2} = y z^{[H-E_1-E_3]}+\frac{z^{[E_1-E_5]}}{x y}+\frac{1}{x},\\
&\vartheta_{D_3} =x y^2 z^{[2 H -2 E_1 - E_2-E_3]}+y z^{[H-E_1-E_2]}+x y z^{[H-E_1]}\\
&\vartheta_{D_4} =x y z^{[2H-E_1-E_2-E_3-E_4]}+\frac{z^{[E_1]}}{y}+x z^{[H-E_4]},
\end{align*}
as well as
\begin{equation*}
\Omega = \frac{dx}{x}\wedge\frac{dy}{y}.
\end{equation*}
We consider K\"ahler classes 
\begin{align*}
&[\omega^{\delta}] = [(1+\delta)H - \frac{1}{2} E_2 - \frac{1}{2}E_5  - \delta \sum_{i=1,3,4} E_i]\\
&=\frac{1}{2} D_2 + \frac{1}{2}D_4 + (1-\delta)D_1 +\delta D_3 
\end{align*}
for sufficiently small $\delta > 0$. In the $\delta \to 0$ limit we have 
\begin{align*}
&[\omega^0_k] = k[H - \frac{1}{2} E_2 - \frac{1}{2}E_5] = k[D_1 + \frac{1}{2} D_2 + \frac{1}{2}D_4]
\end{align*}
and, on $\mathcal{U}$, 
\begin{align*}
 &W^0_k = e^{-3 \pi  k} x y^2+\frac{e^{\pi  k}}{x y^2}+e^{-3 \pi  k} x y+e^{-2 \pi  k} x y+\frac{e^{\pi  k}}{x y}+e^{-2 \pi 
   k} x\\
   &+e^{-2 \pi  k} y+e^{-\pi k} y+\frac{e^{-\pi k}}{y}+\frac{1}{x y}+\frac{1}{x}+\frac{1}{y}. 
\end{align*}
A direct computation shows that we have
\begin{align*}
&\del_x W^0_k = \frac{e^{-3 \pi  k} (y+1) \left(e^{\pi  k}+y\right) \left(x^2 y^2-e^{3 \pi  k}\right)}{x^2 y^2},\\
&\del_y W^0_k = -\frac{e^{-3 \pi  k} \left(e^{2 \pi  k}-x y^2\right) \left(e^{\pi  k} x y+e^{\pi  k} y+e^{2 \pi  k} y+2 e^{2 \pi  k}+2 x y^2+x
   y\right)}{x y^3}
\end{align*}
from which the critical points can be computed exactly,
\begin{align*}
&p_1 = (1, -e^{\pi  k}),\,
p_2 =  \big( -e^{\pi  k}, -\frac{e^{\frac{\pi  k}{2}} \left(e^{\pi 
   k}-1\right)}{\sqrt{-2 e^{\pi  k}+e^{2 \pi  k}+1}}\big),\,
p_3 =   \big( -e^{\pi  k}, \frac{e^{\frac{\pi  k}{2}}, 
   \left(e^{\pi  k}-1\right)}{\sqrt{-2 e^{\pi  k}+e^{2 \pi  k}+1}}\big),\\
& p_4 =  ( e^{\pi  k},  -1),\, 
p_5 = (  e^{\pi  k}, -e^{\frac{\pi  k}{2}}),\,
p_6 =   ( e^{\pi  k},  e^{\frac{\pi 
   k}{2}}),\,p_7= (e^{\pi  k}, -e^{\pi  k}),\,p_8= ( e^{2 \pi  k}, -1).  
\end{align*}
By computations similar to (but much simpler than) \cite{Barrott_equations}, Example 6.1, one can check that indeed all the critical points of $W^0_k$ are contained in $\mathcal{U}$.
\section{Slope stability}\label{SlopeSec}
The Donaldson-Futaki invariant of degeneration to the normal cone $\mathcal{X}$ of a subvariety $Z\subset X$ has been studied in detail in \cite{RossThomas_obstruction}, where it is shown to be equivalent to a quantity called \emph{quotient slope} $\mu_c(Z)$ (recalled below in the surface case), where $c$ corresponds to the additional K\"ahler parameter on $\mathcal{X}$. In particular, $(X, [\omega])$ is (strictly) destabilised by $Z$ precisely when we have
\begin{equation*}
\mu_c(\mathcal{O}_Z) < \mu(X). 
\end{equation*}
\begin{thm}[\cite{RossThomas_obstruction}, Corollary 5.3] Let $Z$ be a smooth curve in a smooth K\"ahler surface $(X,[\omega])$. Then
\begin{align*}
\mu(X) = -\frac{K_X . [\omega]}{[\omega]^2},\, \mu_c(\mathcal{O}_Z) = \frac{3(2[\omega] . Z - c(K_X . Z + Z^2))}{2c(3[\omega] . Z-c Z^2)}.
\end{align*}
\end{thm}
Suppose now $(X, D),\,Z \subset D$ are of the type discussed in the previous Sections. As a special case of the results discussed in Section \ref{DFAsymp}, we see that \emph{if the quantities
\begin{equation*}
\sum_{p \in \operatorname{Crit}(W_k)}\frac{3\vartheta_{Z,k}(2\vartheta_{[\omega],k}  + c(W_k  - \vartheta_{Z,k} ))}{(x y)^2 \det \nabla^2 W_k }\big|_p,\,\sum_{p \in \operatorname{Crit}(W_k)}\frac{2c \vartheta_{Z, k}(3\vartheta_{[\omega],k}   -c \vartheta_{Z,k} )}{(x y)^2 \det \nabla^2 W_k }\big|_p
\end{equation*}
concentrate at a set of critical points $\{p_k\}$, each carrying the same contribution, then we have  
\begin{align*}
&\mu_{c }(\mathcal{O}_Z) = \frac{1}{c}\frac{(1 + 3 c(\frac{k W_k }{\vartheta_{[\omega_k]}} - \frac{k \vartheta_Z }{\vartheta_{[\omega_k]} }))}{ (1  -2 c \frac{k\vartheta_Z }{\vartheta_{[\omega_k]} })}\big|_{\hat{p}_k} + O(k^{-1}), 
\end{align*}
for each such critical point $\hat{p}_k$}.
\subsection{$X = \operatorname{Bl}_p \PP^2$}\label{BlpSlopeSubsec}
Consider $X = \operatorname{Bl}_p \PP^2$ with K\"ahler class 
\begin{align*}
&\omega_k = k (H- q E). 
\end{align*}
The slope is
\begin{equation*}
\mu(X,\omega_k) = \frac{1}{k} \frac{(3H-E).(H-q E)}{(H-q E)^2} = \frac{1}{k} \frac{3-q}{1-q^2}. 
\end{equation*}
We allow an arbitrary scaling $(x, y) \mapsto (e^{-\lambda} x, e^{-\lambda y})$ preserving the symmetry exchanging $x$, $y$. The corresponding Landau-Ginzburg potential is given by
\begin{align*}
&W_k =  \frac{e^{-2\pi k(1-\lambda)}}{x y} +  e^{-\pi k\lambda} x +  e^{-\pi k\lambda} y + e^{2\pi k (q-\lambda) } x y =\vartheta_H + \vartheta_{L_1} + \vartheta_{L_2} + \vartheta_{E}.
\end{align*}
Similarly, we have
\begin{align*}
&\vartheta_{\omega_k} = k(\vartheta_H  - q \vartheta_{E}).  
\end{align*}
The critical points of $W_k$ are solutions in $(\C^*)^2$ of
\begin{align*}
&\del_x W_k =  e^{-\pi  k \lambda }+y e^{2 \pi  k (q-\lambda )}-\frac{e^{-2 \pi  k (1-\lambda )}}{x^2
   y} = 0,\\
&\del_y W_k = e^{-\pi k \lambda }+x e^{2 \pi  k (q-\lambda )}-\frac{e^{-2 \pi  k (1-\lambda )}}{x
   y^2} = 0
\end{align*}
and so they are given by points of $(\C^*)^2$ satisfying
\begin{equation*}
x =y,\, -e^{\pi  k \lambda }+x^4 e^{2 \pi  k (1-\lambda )+\pi  k \lambda +2 \pi  k (q-\lambda
   )}+x^3 e^{2 \pi  k (1-\lambda )} = 0.
\end{equation*}
We claim we can choose the scaling parameter $\lambda$ and the K\"ahler parameter $q$ so that the critical points satisfy either $x \to -1$ (for a single distinguished critical point) or $x \to 0$. The value of $\lambda$ is determined by the condition on the exponents
\begin{equation*}
2 \pi  k (1-\lambda )+\pi  k \lambda +2 \pi  k (q-\lambda) = 2 \pi  k (1-\lambda ) \iff \lambda = 2 q,
\end{equation*}
for which the critical points equation becomes
\begin{equation*}
x^3 (x+1)  = e^{2 \pi  k (3 q - 1)}.
\end{equation*}
Then, the required asymptotic behaviour holds iff $q < \frac{1}{3}$. From now on we fix $\lambda = 2q$ and $q < \frac{1}{3}$. We compute
\begin{align*}
(\det \nabla^2 W_k)^{-1} = -\frac{1}{-\frac{3 e^{4 \pi  k (2 q-1)}}{x^4 y^4}+\frac{2 e^{2 \pi  k (q-1)}}{x^2
   y^2}+e^{-4 \pi  k q}},
\end{align*}
so, at a critical point $(x, x)$, we have
\begin{align*}
(\det \nabla^2 W_k)^{-1} = -\frac{1}{-\frac{3 e^{4 \pi  k (2 q-1)}}{x^8}+\frac{2 e^{2 \pi  k (q-1)}}{x^4}+e^{-4 \pi
    k q}}.
\end{align*}
Now choose $Z = E$ the exceptional divisor. Then we have
\begin{equation*}
\vartheta_Z = \vartheta_E = e^{-2\pi k q} x y
\end{equation*}
and, at a critical point $p = (x,x)$, setting $c = k s$,
\begin{align*}
\frac{\vartheta_Z }{(x y)^2 \det \nabla^2 W }3(2\vartheta_{[\omega]}  + c(W  - \vartheta_Z ))\big|_{p} = \frac{6 k x^7 (q x-s)-3 k (s+2) x^4 e^{2 \pi  k (3 q-1)}}{2 x^4 e^{2 \pi  k (3 q-1)}-3
   e^{4 \pi  k (3 q-1)}+x^8},
\end{align*} 
as well as
\begin{align*}
\frac{\vartheta_Z }{(x y)^2\det \nabla^2 W }2c(3\vartheta_{[\omega]}   -c \vartheta_Z )\big|_{p}= \frac{2 k^2 s x^4 \left(x^4 e^{2 \pi  k (1-3 q)} (3 q+s)-3\right)}{x^8 e^{2 \pi  k (1-3
   q)}-3 e^{2 \pi  k (3 q-1)}+2 x^4}.
\end{align*}
We need to consider the asymptotics of these two terms as $k\to\infty$. In this limit, at a critical point, we have either $x \to 0$ or $x\to -1$. Clearly, the critical points for which $x\to 0$ give a subleading contribution. On the other hand, there is a single critical point for which $x\to -1$, for which we have
\begin{align*}
\frac{\vartheta_Z }{(x y)^2 \det \nabla^2 W }3(2\vartheta_{[\omega]} + c(W  - \vartheta_Z ))\big|_p\sim \frac{3 k \left(2 e^{2 \pi  k (1-3 q)} (q+s)-s-2\right)}{e^{2 \pi  k (1-3 q)}-3 e^{2 \pi
    k (3 q-1)}+2},
\end{align*}
respectively
\begin{align*}
\frac{\vartheta_Z }{(x y)^2 \det \nabla^2 W }2c(3\vartheta_{[\omega]} - c \vartheta_Z )\big|_p\sim \frac{2 k^2 s \left(e^{2 \pi  k (1-3 q)} (3 q+s)-3\right)}{e^{2 \pi  k (1-3 q)}-3 e^{2
   \pi  k (3 q-1)}+2}.
\end{align*}
So, as $k \to \infty$, the leading contribution to the quotient slope of $Z$ is concentrated at the critical point for which $x \to -1$, and is given by
\begin{align*}
\mu_{c_k}(\mathcal{O}_Z) \sim \frac{3 k \left(2 e^{2 \pi  k (1-3 q)} (q+s)-s-2\right)}{2 k^2 s \left(e^{2 \pi  k (1-3 q)} (2 q+s)-2\right)} \sim \frac{1}{k} \frac{3(q+s)}{s(3q + s)}.
\end{align*}
We should compare this to the ambient slope $\mu(X, \omega_k)$. At the threshold value $s = 1-q$, we find
\begin{equation*}
\lim_{k\to \infty}k\mu_{c_k}(\mathcal{O}_Z, \omega_k)\big|_{s=1-q} = \frac{3}{(1-q)(2q + 1)}   < \frac{3-q}{1-q^2} = k\mu(X, \omega_k).
\end{equation*}

Consider now the case when $Z = H$, a line in $\PP^2$. Then, 
\begin{equation*}
\vartheta_Z = \frac{e^{-2 \pi  k (1-2 q)}}{x y}, 
\end{equation*} 
and, at a critical point $p = (x,x)$, we have
\begin{align*}
\frac{\vartheta_Z }{(x y)^2\det \nabla^2 W }3(2\vartheta_{[\omega]} + c(W - \vartheta_Z ))\big|_p = \frac{3 k \left(2 e^{2 \pi  k (3 q-1)}+x^3 (s (x+2)-2 q x)\right)}{x^8 \left(-e^{2 \pi 
   k (1-3 q)}\right)+3 e^{2 \pi  k (3 q-1)}-2 x^4},
\end{align*} 
as well as
\begin{align*}
\frac{\vartheta_Z }{(x y)^2 \det \nabla^2 W}2c(3\vartheta_{[\omega]} -  \vartheta_Z )\big|_p= \frac{2 k^2 s \left(-(s-3) e^{2 \pi  k (3 q-1)}-3 q x^4\right)}{x^8 \left(-e^{2 \pi  k
   (1-3 q)}\right)+3 e^{2 \pi  k (3 q-1)}-2 x^4}.
\end{align*}
The critical points equation shows that we have either $x \to -1$ or $x \sim \xi e^{\frac{2}{3}(3q-1)}$, where $\xi^3 = 1$. Using the explicit formulae we see that the leading contributions to the quantities above is given by critical points for which $x\to 0$, namely
\begin{align*}
\frac{\vartheta_Z }{(x y)^2 \det \nabla^2 W }3(2\vartheta_{[\omega]} + c_k(W - \vartheta_Z ))\big|_p \sim \frac{3 k \left(2 (s+1) e^{4 \pi  k q}+e^{\frac{2}{3} \pi  k (9 q-1)} (s-2 q)\right)}{3
   e^{4 \pi  k q}-e^{\frac{4}{3} \pi  k (6 q-1)}-2 e^{\frac{2}{3} \pi  k (9 q-1)}},
\end{align*}
respectively
\begin{align*}
&\frac{\vartheta_Z }{(x y)^2 \det \nabla^2 W }2c(3\vartheta_{[\omega]} - c \vartheta_Z )\big|_p \sim -\frac{2 e^{\frac{2 \pi  k}{3}} k^2 s \left(3 q e^{2 \pi  k q}+e^{\frac{2 \pi  k}{3}}
   (s-3)\right)}{-e^{4 \pi  k q}-2 e^{\frac{2}{3} \pi  k (3 q+1)}+3 e^{\frac{4 \pi 
   k}{3}}}.
\end{align*}
So we find for the quotient slope of $H$ 
\begin{align*}
& \mu_{c_k}(\mathcal{O}_Z,\omega_k) \sim -\frac{1}{k}\frac{3 \left(e^{2 \pi  k q} (s-2 q)+2 e^{\frac{2 \pi  k}{3}} (s+1)\right)}{2 s \left(3
   q e^{2 \pi  k q}+e^{\frac{2 \pi  k}{3}} (s-3)\right)} \sim \frac{3(s+1)}{s(3-s)}.
\end{align*}
This matches the intersection-theoretic formula, since $\omega = H - q E$, $Z = H$. Note that each critical point contributes equally to the asymptotics (i.e. the latter does not depend on $\xi$), so the assumptions of Theorem \ref{SlopeThm} are satisfied in this example too. 
\subsection{$X = \operatorname{Bl}_{p,q}\PP^2$}\label{BlpqSlopeSubsec}
Let $E$, $F$ denote the exceptional divisors. The toric boundary is given by
\begin{equation*}
D_1 \sim H-F, D_2 \sim H-E, D_3 \sim E, D_4 \sim H-E-F, D_5\sim F.
\end{equation*}
Consider a polarisation
\begin{equation*}
\omega = H - a E - a F = \frac{1}{3} D_1 + \frac{1}{3} D_2 + \big(\frac{2}{3}-a\big)D_3+ \frac{1}{3}D_4+\big(\frac{2}{3}-a\big)D_5. 
\end{equation*}
Thus, 
\begin{align*}
&W_k = \vartheta_{D_1} +  \vartheta_{D_2} +  \vartheta_{D_3} +  \vartheta_{D_4}+  \vartheta_{D_5}\\
&= x + y + e^{2\pi k (a-1)} \frac{1}{x} + e^{-2\pi k} \frac{1}{x y} + e^{2\pi k (a-1)} \frac{1}{y}.
\end{align*}
Similarly, 
\begin{equation*}
\vartheta_{[\omega_k]} = \frac{k}{3}\vartheta_{D_1} + \frac{k}{3} \vartheta_{D_2} + \big(\frac{2}{3}-a\big)k\vartheta_{D_3}+ \frac{k}{3}\vartheta_{D_4}+\big(\frac{2}{3}-a\big)k\vartheta_{D_5}.
\end{equation*}
Now choose 
\begin{equation*}
a = \frac{1}{2} - \delta
\end{equation*}
for sufficiently small $\delta > 0$. We have
\begin{align*}
&\del_x W_k =  -\frac{e^{-\pi k (2 \delta +1) }}{x^2}-\frac{e^{-2 \pi  k}}{x^2 y}+1,\\
&\del_y W_k = -\frac{e^{-\pi k  (2 \delta +1)}}{y^2} -\frac{e^{-2 \pi  k}}{x y^2}+1 ,
\end{align*}
and we find that critical points are solutions of
\begin{align*}
& \left(e^{4 \pi  \left(-\delta -\frac{1}{2}\right) k+2 \pi  k}+x^2 \left(-e^{2 \pi  \left(-\delta -\frac{1}{2}\right) k+2 \pi  k}\right)-x\right)
   \left(-e^{2 \pi  k} x^3+x e^{2 \pi  \left(-\delta -\frac{1}{2}\right) k+2 \pi  k}+1\right) = 0,\\
&y = -\frac{e^{-2 \pi  k}}{e^{2 \pi  \left(-\delta -\frac{1}{2}\right) k}-x^2}.
\end{align*}
From this we can read off the asymptotics of the critical points
\begin{align*}
& p_1 \sim (-e^{\pi k (-1+2 \delta)} ,  -e^{\pi k (-1+2 \delta)}),\\
& p_i \sim (\pm e^{\pi k (-\frac{1}{2}-\delta)}, e^{\pi k (-\frac{1}{2}-\delta)} ),\,i = 2, 3,\, p_i \sim (\pm e^{\pi k (-\frac{1}{2}-\delta)} , -e^{\pi k (-\frac{1}{2}-\delta)}),\,i = 4, 5.
\end{align*}
Choose $Z = D_4 \sim H-E-F$, so $\vartheta_Z = \frac{e^{-2\pi k}}{x y}$. Set
\begin{align*}
& d(p) = \frac{\vartheta_Z }{(x y)^2 \det \nabla^2 W }3(2\vartheta_{[\omega]} + c(W - \vartheta_Z ))|_p,\\
& r(p) = \frac{\vartheta_Z }{(x y)^2 \det \nabla^2 W }2c(3\vartheta_{[\omega]} - c \vartheta_Z )|_p.
\end{align*}
Using the above critical asymptotics, we compute
\begin{align*}
& d(p_1) = 2 k \left((3 s+2) e^{\pi  (6 \delta -1) k}+3 (2 \delta +s)\right) \sim 6 k  (s + 2 \delta ),\\
& d(p_2) = \frac{2 k \left(e^{\frac{1}{2} \pi  (1-6 \delta ) k} (6 \delta +6 s+3)+1\right)}{8 e^{\frac{1}{2} \pi  (1-6 \delta ) k}+4 e^{\pi  (1-6 \delta )
   k}+3} \to 0,\\
& d(p_3) = -\frac{2 k}{4 e^{\pi  (1-6 \delta ) k}-3} \to 0,\\
& d(p_4) = -\frac{2 k}{4 e^{\pi  (1-6 \delta ) k}-3} \to 0,\\
& d(p_5) = -\frac{2 k \left(e^{\frac{1}{2} \pi  (1-6 \delta ) k} (6 \delta +6 s+3)-1\right)}{-8 e^{\frac{1}{2} \pi  (1-6 \delta ) k}+4 e^{\pi  (1-6 \delta )
   k}+3} \to 0,
\end{align*} 
respectively
\begin{align*}
& r(p_1) = 2 k^2 s \left(6 \delta +2 e^{\pi  (6 \delta -1) k}+s\right) \sim 2 k^2 s(s + 6 \delta),\\
& r(p_2) = \frac{2 k^2 s \left((6 \delta +3) e^{\frac{1}{2} \pi  (1-6 \delta ) k}-s+1\right)}{8 e^{\frac{1}{2} \pi  (1-6 \delta ) k}+4 e^{\pi  (1-6 \delta )
   k}+3} \to 0,\\
& r(p_3) = \frac{2 k^2 (s-1) s}{4 e^{\pi  (1-6 \delta ) k}-3} \to 0,\\   
& r(p_4) = \frac{2 k^2 (s-1) s}{4 e^{\pi  (1-6 \delta ) k}-3} \to 0,\\
& r(p_5) =-\frac{2 k^2 s \left((6 \delta +3) e^{\frac{1}{2} \pi  (1-6 \delta ) k}+s-1\right)}{-8 e^{\frac{1}{2} \pi  (1-6 \delta ) k}+4 e^{\pi  (1-6 \delta )
   k}+3} \to 0.
\end{align*}
It follows that $p_1$ is the single critical point yielding the leading contribution to the quotient slope as $k \to \infty$, given by
\begin{align*}
&\mu_{k s}(\mathcal{O}_Z, \omega_k) \sim \frac{1}{k} \frac{3(s + 2\delta)}{s(s+6\delta)}.
\end{align*}

\subsection{Simple del Pezzo degeneration}\label{IteratedBlpSlopeSec}

Let $X$ be the toric manifold obtained from $\operatorname{Bl}_p \PP^2$ by blowing up a torus fixed point on the exceptional divisor. $X$ contains distinguished rational curves $E_1, E_2$ with
\begin{equation*}
E^2_1 = -2,\, E^2_2 = -1,\,E_1 . E_2 = 1. 
\end{equation*}
Then we have $D = \sum^5_{i = 1} D_i$ with
\begin{equation*}
D_1 \sim H, D_2 \sim H - E_1 - 2 E_2,\,D_3 \sim E_2,\,D_4 \sim E_1,\,D_5 \sim H - E_1 - E_2. 
\end{equation*}
Note that this is weak del Pezzo: the anticanonical bundle $-K_X \sim 3H - E_1 - 2 E_2$ is nef. Indeed, $X$ is a resolution of the simple degeneration $X'$ of the blowup $\operatorname{Bl}_{p, q} \PP^2$ obtained when $p, q$ collide. The class
\begin{equation*}
[\omega] = H - \frac{1}{2} E_1 - r E_2 = D_1 - r D_3 - \frac{1}{2} D_4    
\end{equation*}
is K\"ahler for $1 > r > \frac{1}{2}$, and so we have, with respect to the rescaled class $k [\omega]$,
\begin{align*}
\vartheta_{D_1} = x,\,\vartheta_{D_2} =  \frac{e^{-2 \pi k}}{x y} ,\,\vartheta_{D_3} = \frac{e^{2 k \pi (-1 + r)}}{x^2 y},\,\vartheta_{D_4} = \frac{e^{-k \pi}}{x},\,\vartheta_{D_5} = y ,
\end{align*} 
respectively
\begin{align*}
&W_k =  x + \frac{e^{-2 \pi k}}{x y} + \frac{e^{2 k \pi (-1 + r)}}{x^2 y} + \frac{e^{-k \pi}}{x} + y,\\
&\vartheta_{[\omega_k]} = x   -\frac{1}{2} \frac{e^{-k \pi}}{x} - r \frac{e^{2 k \pi (-1 + r)}}{x^2 y}.
\end{align*}
The critical locus is given by
\begin{align*}
&\del_x W = -\frac{2 e^{2 \pi  k (r-1)}}{x^3 y}-\frac{e^{-2 \pi  k}}{x^2 y}-\frac{e^{\pi  (-k)}}{x^2}+1,\\
& \del_yW = -\frac{e^{2 \pi  k (r-1)}}{x^2 y^2}-\frac{e^{-2 \pi  k}}{x y^2}+1.
\end{align*}
Along the expansion for critical points 
\begin{equation*}
x = a e^{ \pi k b}(1+ O(k^{-1})),\,y = c e^{ \pi k d}(1+O(k^{-1})),
\end{equation*}
we compute
\begin{align*}
&\del_x W \sim -\frac{2 e^{-3 \pi  b k-\pi  d k+2 \pi  k (r-1)}}{a^3 c}-\frac{e^{-2 \pi  b k-\pi  d k-2 \pi  k}}{a^2 c}-\frac{e^{-2 \pi  b k-\pi  k}}{a^2}+1,\\
& \del_yW \sim -\frac{e^{-2 \pi  b k-2 \pi  d k+2 \pi  k (r-1)}}{a^2 c^2}-\frac{e^{-\pi  b k-2 \pi  d k-2 \pi  k}}{a c^2}+1.
\end{align*}
Let us consider the potential asymptotic behaviour of critical points determined by 
\begin{align*}
&\{-3 \pi  b k-\pi  d k+2 \pi  k (r-1)=0,-2 \pi  b k-2 \pi  d k+2 \pi  k (r-1)=0\} \\
&\iff  b = \frac{r-1}{2},\,d = \frac{r-1}{2}.
\end{align*}
Then, in this chamber, we have
\begin{align*}
&\del_x W \sim -\frac{2}{a^3 c}-\frac{e^{-\frac{3}{2} \pi  k (r-1)-2 \pi  k}}{a^2 c}-\frac{e^{-\pi  k (r-1)-\pi  k}}{a^2}+1,\\
& \del_yW \sim -\frac{1}{a^2 c^2}-\frac{e^{-\frac{3}{2} \pi  k (r-1)-2 \pi  k}}{a c^2}+1.
\end{align*}
So we see that this asymptotic behaviour is realised for all $r \in (\frac{1}{2},1)$, with coefficients determined by
\begin{equation*}
-\frac{2}{a^3 c} +1 = -\frac{1}{a^2 c^2} +1 = 0 \iff (a, c) \in \{(\pm 2^{1/2}, \pm 2^{-1/2}),\,(\pm 2^{1/2} \ii, \pm 2^{-1/2}\ii)\}. 
\end{equation*}
This give a set of critical points 
\begin{align*}
& p_i \sim (\pm 2^{1/2} e^{-\frac{\pi k}{2}(1-r)}, \pm 2^{-1/2}e^{-\frac{\pi k}{2}(1-r)}),\,i = 2, 3,\\
& p_i \sim (\pm 2^{1/2} \ii e^{-\frac{\pi k}{2}(1-r)}, \pm 2^{-1/2}\ii e^{-\frac{\pi k}{2}(1-r)}),\,i = 4, 5.
\end{align*}
On the other hand, solving $\del_x W = 0$ with respect to $y$, we find
\begin{equation*}
y = \frac{e^{\pi  (-k)} \left(2 e^{2 \pi  k (r-1)+2 \pi  k}+x\right)}{x \left(e^{\pi  k} x^2-1\right)}
\end{equation*} 
and substituting into $\del_y W$ shows that the critical points equation is equivalent to the quintic
\begin{equation*}
- x \left(e^{\pi  k} x^2-1\right)^2- e^{2 \pi  k r} \left(e^{\pi  k} x^2-1\right)^2 +\left(2 e^{2 \pi  k r}+x\right)^2 = 0,
\end{equation*}
which has three distinct real roots. Thus, there is an additional distinguished critical point $p_1$. 

Note that we have 
\begin{equation*}
\int_{E_1} \omega = 1 - r \to 0 \text{ as } r \to 1 
\end{equation*}
so the limiting class for $r = 1$ contracts $E_1$, going back to the simple degeneration $X'$. According to \cite{RossThomas_obstruction}, Example 5.34, the $-2$ curve $E_1$ slope-destabilises $X$ for $r$ close to $1$, for a suitable choice of parameter $c_r$. Let us show that the assumptions of Theorem \ref{SlopeThm} hold in this case.

So choose $Z = D_4 \sim E_1$, $\vartheta_Z = \frac{e^{-k \pi}}{x}$. Set
\begin{align*}
& d(p) = \frac{\vartheta_Z }{(x y)^2 \det \nabla^2 W }3(2\vartheta_{[\omega]} + c(W - \vartheta_Z ))|_{p},\\
& r(p) = \frac{\vartheta_Z }{(x y)^2 \det \nabla^2 W }2c(3\vartheta_{[\omega]} - c \vartheta_Z )|_{p}.
\end{align*}
Using the above critical asymptotics, we compute in the large $k$ limit  
\begin{align*}
& d(p_2)|_{r=1} = \frac{3 e^{\pi  k} k \left(e^{2 \pi  k} (4 s+2)-e^{\pi  k}+\sqrt{2} s\right)}{4 \sqrt{2} e^{\pi  k}+12 \sqrt{2} e^{2 \pi  k}+4 e^{3 \pi  k}+8 e^{4 \pi  k}+6} \to 0,\\
& d(p_3)|_{r=1} = \frac{3 e^{\pi  k} k \left(e^{2 \pi  k} (4 s+2)-e^{\pi  k}-\sqrt{2} s\right)}{-4 \sqrt{2} e^{\pi  k}-12 \sqrt{2} e^{2 \pi  k}+4 e^{3 \pi  k}+8 e^{4 \pi  k}+6} \to 0,\\
& d(p_4)|_{r=1} = -\frac{3 e^{\pi  k} k \left(e^{2 \pi  k} (4 s+2)+e^{\pi  k}+\ii \sqrt{2} s\right)}{-4 i \sqrt{2} e^{\pi  k}+12 \ii \sqrt{2} e^{2 \pi  k}-4 e^{3 \pi  k}+8 e^{4 \pi  k}-6} \to 0,\\
& d(p_5)|_{r=1} = -\frac{3 e^{\pi  k} k \left(e^{2 \pi  k} (4 s+2)+e^{\pi  k}-\ii \sqrt{2} s\right)}{4 \ii \sqrt{2} e^{\pi  k}-12 \ii \sqrt{2} e^{2 \pi  k}-4 e^{3 \pi  k}+8 e^{4 \pi  k}-6} \to 0,
\end{align*} 
respectively
\begin{align*}
& r(p_2)|_{r=1} = \frac{e^{2 \pi  k} k^2 s \left(6 e^{\pi  k}-2 s-3\right)}{4 \sqrt{2} e^{\pi  k}+12 \sqrt{2} e^{2 \pi  k}+4 e^{3 \pi  k}+8 e^{4 \pi  k}+6} \to 0,\\
& r(p_3)|_{r=1} = \frac{e^{2 \pi  k} k^2 s \left(6 e^{\pi  k}-2 s-3\right)}{-4 \sqrt{2} e^{\pi  k}-12 \sqrt{2} e^{2 \pi  k}+4 e^{3 \pi  k}+8 e^{4 \pi  k}+6} \to 0,\\   
& r(p_4)|_{r=1} = \frac{e^{2 \pi  k} k^2 s \left(6 e^{\pi  k}+2 s+3\right)}{4 \ii \sqrt{2} e^{\pi  k}-12 \ii \sqrt{2} e^{2 \pi  k}+4 e^{3 \pi  k}-8 e^{4 \pi  k}+6} \to 0,\\
& r(p_5)|_{r=1} = \frac{e^{2 \pi  k} k^2 s \left(6 e^{\pi  k}+2 s+3\right)}{-4 \ii \sqrt{2} e^{\pi  k}+12 \ii \sqrt{2} e^{2 \pi  k}+4 e^{3 \pi  k}-8 e^{4 \pi  k}+6} \to 0.
\end{align*} 
Thus, for $r$ sufficiently close to $1$, we must have
\begin{equation*}
\mu_{k s}(\mathcal{O}_Z, \omega_k) \sim \frac{1}{k}\frac{d(p_1)}{r(p_1)}.
\end{equation*}
\subsection{Degree $5$ del Pezzo}\label{Deg5SlopeSec}
Continuing the non-toric example of Section \ref{Deg5AsympSec}, we compute in this case
\begin{align*}
\frac{1}{(x y)^2 \det \nabla^2 W^0} = \frac{e^{8\pi k} x^2 y^2}{Q},
\end{align*}
where
\begin{align*}
& Q = 3 -e^{6 \pi  k} x \left(x^3-4 y^2\right)+2 e^{7 \pi  k} x^2 y^2+4 e^{5 \pi  k} x y (2
   x+y)+2 e^{4 \pi  k} y (4 x+3 y)\\&+4 e^{2 \pi  k} (x+2 y)+2 e^{3 \pi  k} x (3 x+4 y)+4
   e^{\pi  k} x-e^{8 \pi  k} y^4.
\end{align*}
Consider the boundary component $Z = D_1$ with corresponding theta function $\vartheta_{Z} = \vartheta_1$. Note that we have
\begin{align*}
& [\omega^0] = H - \frac{1}{2} E_1 = \frac{1}{2} D_1 + \frac{1}{2} D_2 + \frac{1}{2} D_3 + \frac{1}{2} D_4,\\ 
& \vartheta_{[\omega^0]} = \frac{1}{2} \vartheta_1 + \frac{1}{2} \vartheta_2 + \frac{1}{2} \vartheta_3 + \frac{1}{2} \vartheta_4\\
&= \frac{1}{2} x + \frac{1}{2} y + \frac{1}{2} \big(\frac{e^{-2 \pi  k}}{x}+\frac{y}{x}\big) + \frac{1}{2} \big(\frac{e^{-4 \pi  k}}{x y}+\frac{e^{-2 \pi  k}}{x}+\frac{e^{-3 \pi  k}}{y}\big).
\end{align*}
Then, for $\delta = 0$, and at points $p\in\operatorname{Crit}(W^0) \cap \mathcal{U}$, we have 
\begin{align*}
&\frac{\vartheta_Z }{( x y )^2 \det \nabla^2 W^0 }3(2k\vartheta_{[\omega^0]} + c_k(W^0 - \vartheta_Z ))|_p = \frac{e^{8\pi k} x^3 y^2}{Q} 3(2 k \vartheta_{[\omega^0]} + c_k(W^0 - x))|_p \to 0
\end{align*} 
except for the single critical point with asymptotics $x \sim -1, y\sim 1$. Similarly, we have
\begin{align*}
&\frac{\vartheta_Z }{(x  y)^2 \det \nabla^2 W }2c(3k\vartheta_{[\omega^0]} - c_k \vartheta_Z )|_p= \frac{e^{8\pi k} x^3 y^2}{Q} 2 c_k(3 k \vartheta_{[\omega^0]} -c _k x )|_p \to 0
\end{align*}
except for $x \sim -1, y\sim 1$. Thus, the contributions of critical points to $\mu_c(\mathcal{O}_Z)$ concentrate at $x \sim -1, y\sim 1$, at least in the chart $\mathcal{U}$. 

It remains to estimate contributions from the other charts. However, for these other critical points $p \notin \mathcal{U}$ we have either $\vartheta_{1}(p) = 0$ (which obviously does not contribute), or $\vartheta_{2}(p) = 0$, and we note that the above rational functions of $x, y$ extend continuously to $0$ across the locus $y = \vartheta_{2} = 0$. 
\subsection{Degree $4$ del Pezzo}\label{Deg4SlopeSec} 
Recall from \ref{Deg4AsympSec} that we are fixing K\"ahler classes $[\omega^{\delta}]$ with limit $[\omega^0] = [D_1 + \frac{1}{2} D_2 + \frac{1}{2}D_4]$, so we have
\begin{equation*}
\vartheta_{[\omega^0_k]} = k\big(\vartheta_{D_1} + \frac{1}{2} D_2 + \frac{1}{2}D_3\big).
\end{equation*}
Choosing $Z = D_1$, with theta function (corresponding to $[\omega^0_k]$)
\begin{equation*}
\vartheta_{D_1} = \frac{e^{\pi k}}{x y^2}+\frac{e^{-\pi k}}{y}+\frac{1}{x y},
\end{equation*}
we set 
\begin{align*}
&\tilde{d}(p) = k^{-1}\frac{\vartheta_{D_1} }{(x y)^2 \det \nabla^2 W^0_k }3(2 k\vartheta_{[\omega^0_k]} + k s(W^0_k - \vartheta_{D_4} ))\big|_p.
\end{align*}
Using the critical asymptotics given in Section \ref{Deg4AsympSec}, we compute 
\begin{align*}
&\tilde{d}(p_1) = -\frac{3   \left(e^{\pi  k} \left(e^{\pi  k} (s+1)+2 s+1\right)+2\right)}{\left(e^{\pi  k}-1\right)^4} \to 0,\\
&\tilde{d}(p_2) = \frac{3 e^{\frac{\pi  k}{2}}   \left(e^{\pi  k}-1\right)   -3 e^{\pi  k}   (3 s+4)}{4 \left(e^{\pi  k}-1\right)^2}\to0,\\
&\tilde{d}(p_3) = -\frac{3 e^{\frac{\pi  k}{2}}   \big(e^{\frac{\pi  k}{2}} (3 s+4)+ \left(e^{\pi  k}-1\right) \big)}{4 \left(e^{\pi 
   k}-1\right)^2} \to 0,\\
&\tilde{d}(p_4) = \frac{3 e^{\pi  k} \left(e^{\pi  k}-2\right) \left(e^{2 \pi  k} (2 s-1)+5 e^{\pi  k}+s\right)}{\left(e^{\pi  k}-1\right)^4} \to 3(2s-1),\\
&\tilde{d}(p_5) = -\frac{3 \left(e^{\frac{\pi  k}{2}}-2\right) \left(-e^{\frac{\pi  k}{2}} (3 s+4)+e^{\pi  k} (4 s+3)+2 s+5\right)}{4 \left(e^{\frac{\pi k}{2}}-1\right)^4} \to 0,\\
\end{align*}
\begin{align*}
&\tilde{d}(p_6) = \frac{3 \left(e^{\frac{\pi  k}{2}}+2\right) \left(e^{\frac{\pi  k}{2}} (3 s+4)+e^{\pi  k} (4 s+3)+2 s+5\right)}{4 \left(e^{\frac{\pi 
   k}{2}}+1\right)^4} \to 0,\\
&\tilde{d}(p_7) = d(p_8) = -\frac{3 \left(e^{\pi  k} \left(e^{\pi  k} (s+1)+2 s+1\right)+2\right)}{\left(e^{\pi  k}-1\right)^4} \to 0.   
\end{align*}
Similarly, setting
\begin{align*}
&\tilde{r}(p) = k^{-2} \frac{\vartheta_{D_1} }{(x y)^2 \det \nabla^2 W }2 k s(3k\vartheta_{[\omega_k]} - k s \vartheta_{D_1} )\big|_{p},
 \end{align*}
 we compute
\begin{align*}
&\tilde{r}(p_1) = -\frac{s \left(3 e^{\pi  k} \left(e^{\pi  k}+1\right)-2 s+6\right)}{\left(e^{\pi  k}-1\right)^4} \to 0,\\
&\tilde{r}(p_2) =\frac{e^{\frac{\pi  k}{2}} s \left(2 e^{\frac{\pi  k}{2}} (s-6)+3  \left(e^{\pi  k}-1\right) \right)}{4 \left(e^{\pi 
   k}-1\right)^2} \to 0,\\
&\tilde{r}(p_3) = \frac{e^{\frac{\pi  k}{2}}  s \left(2 e^{\frac{\pi  k}{2}} (s-6)-3  \left(e^{\pi  k}-1\right) \right)}{4 \left(e^{\pi 
   k}-1\right)^2}  \to 0,\\
&\tilde{r}(p_4) = \frac{e^{2 \pi  k} \left(e^{\pi  k}-2\right) s \left(e^{\pi  k} (2 s-3)-4 s+15\right)}{\left(e^{\pi  k}-1\right)^4} \to s(2s-3),\\
&\tilde{r}(p_5) = -\frac{\left(e^{\frac{\pi  k}{2}}-2\right) s \left(2 e^{\frac{\pi  k}{2}} (s-6)+9 e^{\pi  k}-4 s+15\right)}{4 \left(e^{\frac{\pi k}{2}}-1\right)^4}  \to 0,\\
&\tilde{r}(p_6) = \frac{\left(e^{\frac{\pi  k}{2}}+2\right) s \left(-2 e^{\frac{\pi  k}{2}} (s-6)+9 e^{\pi  k}-4 s+15\right)}{4 \left(e^{\frac{\pi k}{2}}+1\right)^4} \to 0,\\
&\tilde{r}(p_7) = \tilde{r}(p_8) = -\frac{ s \left(3 e^{\pi  k} \left(e^{\pi  k}+1\right)-2 s+6\right)}{\left(e^{\pi  k}-1\right)^4} \to 0.
\end{align*}
It follows that we have
\begin{equation*}
\mu_{ s}(\mathcal{O}_{D_1}, \omega^0) = \lim_{k \to \infty} \frac{\tilde{d}(p_4)}{\tilde{r}(p_4)}   
\end{equation*}
so Theorem \ref{SlopeThm} also applies in this case.
\section{Localised Futaki character}\label{ProductTestSec}
Suppose $V_a$ is a real holomorphic vector field corresponding to an element $a$ lying in a maximal compact torus of autmorphisms of $(X, \omega)$. Following the notation of \cite[Section 3]{Legendre_localised}, we have
\begin{align*}
&\mathcal{F}_{[\omega]}(J V_a) = \sum_Z \int_Z \frac{n \bar{c}_{[\omega]} ([\omega] - \mu_{a}(Z))^{n+1}}{(n+1)! \operatorname{e}(N^Z_a)}\\
&-2\sum_Z \int_Z \frac{(2\pi c_1(X) + \sum^{n-n_Z}_{i=1} \langle \operatorname{w}^Z_i, a\rangle)  \cup ([\omega] - \mu_{a}(Z))^n}{n! \operatorname{e}(N^Z_a)},  
\end{align*}
where we may assume that $Z$ are fixed loci of $V_a$ and $\mu$ is a moment map for a maximal torus containing $V_a$. 

Moreover, if we normalise the moment map by
\begin{equation}\label{momentNormalisation}
\int_X \mu\,\omega^n = 0
\end{equation}
then we have in fact 
\begin{align*}
\mathcal{F}_{[\omega]}(J V_a) = -2\sum_Z \int_Z \frac{(2\pi c_1(X) + \sum^{n-n_Z}_{i=1} \langle \operatorname{w}^Z_i, a\rangle)  \cup ([\omega] - \mu_{a}(Z))^n}{n! \operatorname{e}(N^Z_a)},  
\end{align*}
(see \cite[Remark 3.2]{Legendre_localised}).
We apply this to the case of a complex surface, assuming the normalisation \eqref{momentNormalisation}. Then, the contribution from fixed divisors $D \subset X$ is given by
\begin{align*}
&-2\sum_D \int_D (-2\pi)\frac{(2\pi c_1(X) + \langle \operatorname{w}^D, a\rangle)  \cup ([\omega] - \mu_{a}(D))^2 }{n! (2\pi c_1(N^D) - \langle \operatorname{w}^D, a\rangle)}\\
&= -\frac{4\pi}{n!} \sum_D \frac{1}{\langle \operatorname{w}^D, a\rangle} \int_D (2\pi c_1(X) + \langle \operatorname{w}^D, a\rangle)  \cup ([\omega] - \mu_{a}(D))^2\cup \big(1+\frac{2\pi  c_1(N^D)}{\langle \operatorname{w}^D, a\rangle}\big).
\end{align*}
We compute
\begin{align*}
&\int_D (2\pi c_1(X) + \langle \operatorname{w}^D, a\rangle)  \cup ([\omega] - \mu_{a}(D))^2\cup \big(1+\frac{2\pi  c_1(N^D)}{\langle \operatorname{w}^Z, a\rangle}\big) \\
& =  (\mu_a(D))^2 2\pi c_1(X) . [D] - 2 \langle \operatorname{w}^D, a\rangle \mu_a(D) [\omega] . [D] + (\mu_a(Z))^2 2\pi c_1(N^D) . [D].   
\end{align*}
Now, by adjunction, 
\begin{equation*}
c_1(N^D) . [D] = [D] . [D],
\end{equation*}
so we can write
\begin{align*}
&\int_Z (2\pi c_1(X) + \langle \operatorname{w}^D, a\rangle)  \cup ([\omega] - \mu_{a}(D))^2\cup \big(1+\frac{2\pi  c_1(N^D)}{\langle \operatorname{w}^D, a\rangle}\big) \\
& =  (\mu_a(D))^2 2\pi c_1(X) . [D] - 2 \langle \operatorname{w}^D, a\rangle \mu_a(D) [\omega] . [D] + (\mu_a(D))^2 2\pi [D] . [D],       
\end{align*}
and we find the overall contribution from fixed divisors, 
\begin{align*}
&-4\pi^2 c_1(X) . \big(\sum_D \frac{(\mu_a(D))^2}{\langle \operatorname{w}^D, a\rangle}[D]\big)\\
&+ 4\pi [\omega] . \big(\sum_D \mu_a(D) [D]\big) -4\pi^2 \sum_D \frac{(\mu_a(D))^2}{\langle \operatorname{w}^D, a\rangle} [D].[D].  
\end{align*}
\begin{rmk} This is indeed homogenous of degree $2$ in the K\"ahler moduli $[\omega]$ since $\mu_a(D)$ scales linearly with $[\omega]$.  
\end{rmk}
Suppose now that we have
 \begin{align*}
-2\sum_{\dim Z = 0} \int_Z \frac{(2\pi c_1(X) + \sum^{n-n_Z}_{i=1} \langle \operatorname{w}^Z_i, a\rangle)  \cup ([\omega] - \mu_{a}(Z))^n}{n! \operatorname{e}(N^Z_a)} = 0,  
\end{align*}
and that the only divisor contributions come from fixed divisors $D, D'$. Then, we have
\begin{align*}
&\mathcal{F}_{[\omega]}(J V_a) = -4\pi^2 c_1(X) . \big(\frac{(\mu_a(D))^2}{\langle \operatorname{w}^D, a\rangle}[D] + \frac{(\mu_a(D'))^2}{\langle \operatorname{w}^{D'}, a\rangle}[D'] \big)\\
&+ 4\pi [\omega] . \big( \mu_a(D) [D] + \mu_a(D') [D']\big) -4\pi^2  \frac{(\mu_a(D))^2}{\langle \operatorname{w}^D, a\rangle} [D].[D] -4\pi^2  \frac{(\mu_a(D'))^2}{\langle \operatorname{w}^{D'}, a\rangle} [D'].[D'].  
\end{align*}
So, we can write the vanishing condition $\mathcal{F}_{[\omega]}(J V_a) = 0$ on the mirror as 
\begin{align*}
&\sum_{p\in\operatorname{Crit}(W_k)} \frac{\vartheta_{D,k} \big(-4\pi^2 \frac{(\mu_a(D))^2}{\langle \operatorname{w}^D, a\rangle} W_k + 4\pi \mu_a(D) \vartheta_{[\omega_k]} - 4\pi^2  \frac{(\mu_a(D))^2}{\langle \operatorname{w}^D, a\rangle} \vartheta_{D,k} \big)}{(x y)^2 \det \nabla^2 W_k }\big|_p \\
& = \sum_{p\in\operatorname{Crit}(W_k)} \frac{\vartheta_{D',k} \big(-4\pi^2 \frac{(\mu_a(D'))^2}{\langle \operatorname{w}^{D'}, a\rangle} W_k + 4\pi \mu_a(D') \vartheta_{[\omega_k]} -4\pi^2  \frac{(\mu_a(D'))^2}{\langle \operatorname{w}^{D'}, a\rangle} \vartheta_{D',k} \big)}{(x y)^2 \det \nabla^2 W_k }\big|_p.
 \end{align*} 
 Suppose these quantities concentrate uniformly at sets of critical points $\{p\}, \{p'\}$ in the large $k$ limit. Then, we must have
 \begin{align*}
&\frac{\#\{p\}}{\#\{p'\}} \frac{\mu_a(D)}{\mu_a(D')} \lim_{k\to \infty} \kappa(p,p') \frac{\vartheta_{D,k} (1-\pi \frac{\mu_a(D)}{\langle \operatorname{w}^D, a\rangle} (\frac{W_k }{\vartheta_{[\omega_k]} } + \frac{\vartheta_{D,k} }{\vartheta_{[\omega_k]} })|_p}{\vartheta_{D',k} (1-\pi \frac{\mu_a(D'))}{\langle \operatorname{w}^{D'}, a\rangle} (\frac{W_k }{\vartheta_{[\omega_k]} } + \frac{\vartheta_{D',k} }{\vartheta_{[\omega_k]} }))|_{p'}} = 1,
 \end{align*} 
where
\begin{align*}
&\kappa(p,p') = \frac{\vartheta_{[\omega_k]}(p)}{\vartheta_{[\omega_k]}(p')} \frac{(p'_1 p'_2)^2 \det \nabla^2 W_k(p')}{(p_1 p_2)^2 \det \nabla^2 W_k(p)} \\
&= \frac{\vartheta_{[\omega_k]}(p')}{\vartheta_{[\omega_k]}(p)} \frac{\vartheta^2_{[\omega_k]}(p)}{(p_1 p_2)^2 \det \nabla^2 W_k(p)} \frac{(p'_1 p'_2)^2 \det \nabla^2 W_k(p')}{\vartheta^2_{[\omega_k]}(p')} = \frac{\vartheta_{[\omega_k]}(p')}{\vartheta_{[\omega_k]}(p)} \frac{\vol_{p,k}(\omega)}{\vol_{p', k}(\omega)}.
\end{align*}
More generally, suppose only that the contribution from $D'$ concentrates uniformly on a set of critical points $\{p'\}$. Then, $\mathcal{F}_{[\omega]}(J V_a) = 0$ implies the identity
 \begin{align*}
& \sum_{p\in\operatorname{Crit}(W_k)} \lim_{k\to \infty}\kappa(p,p') \frac{\vartheta_{D,k} (1-\pi \frac{\mu_a(D)}{\langle \operatorname{w}^D, a\rangle} (\frac{W_k }{\vartheta_{[\omega_k]} } + \frac{\vartheta_{D,k} }{\vartheta_{[\omega_k]} }))|_p}{\vartheta_{D',k} (1-\pi \frac{\mu_a(D')}{\langle \operatorname{w}^{D'}, a\rangle} (\frac{W_k }{\vartheta_{[\omega_k]} } + \frac{\vartheta_{D',k} }{\vartheta_{[\omega_k]} }))|_{p'}} =  \#\{p'\} \frac{\mu_a(D')}{\mu_a(D)}.
 \end{align*} 
\subsection{$X = \operatorname{Bl}_p \PP^2$}
Let us consider $X = \operatorname{Bl}_p \PP^2$ where $p = [1:0:0]$ and we choose the holomorphic vector field $V_a$ on $X$ induced by the action on $\PP^2$ given by $t\cdot[z_0:z_1:z_2] = [z_0:t z_1:t z_2]$. The fixed locus consists of $H$, a line in $\PP^2\setminus p$, and $E$, the exceptional divisor, so we are in the situation described above. Continuing the example discussed in Section \ref{BlpSlopeSubsec}, we have, at a critical point $p = (x,x)$,  
\begin{align*}
&\frac{\vartheta_H }{(x y)^2\det \nabla^2 W }4\pi \mu_a(H) \big(\vartheta_{[\omega_k]} - \pi \frac{ \mu_a(H) }{\langle \operatorname{w}^H, a\rangle} W_k - \pi \frac{ \mu_a(H) }{\langle \operatorname{w}^H, a\rangle} \vartheta_{H,k} \big)|_p \\
& = -\frac{4 \pi  \mu_a(H)  \left(x^3 (q \langle \operatorname{w}^H, a\rangle x+\pi  \mu_a(H)  (x+2))-e^{2 \pi  k (3 q-1)} (\langle \operatorname{w}^H, a\rangle-2 \pi  \mu_a(H) )\right)}{\langle \operatorname{w}^H, a\rangle \left(x^8
   \left(-e^{2 \pi  k (1-3 q)}\right)+3 e^{2 \pi  k (3 q-1)}-2 x^4\right)},
\end{align*}
and, as in our previous computations, this concentrates at critical points with asymptotics $x \sim \xi e^{-\frac{1}{3}\pi k}$, where $\xi^3 = 1$. Suppose $q > \frac{1}{9}$. Then the contribution from such a point is
\begin{align*}
&\lim_{k \to \infty}\frac{4 \pi  \mu_a(H)  e^{\frac{2}{3} \pi  k (4-3 q)} \left(2 \pi  e^{\pi  k} \mu_a(H)  \xi ^3+e^{\frac{2 \pi  k}{3}} \xi ^4
   (\pi  \mu_a(H) +q w)-e^{6 \pi  k q} (w-2 \pi  \mu_a(H) )\right)}{w \left(\xi ^8 e^{4 \pi  k (1-2 q)}+2 \xi ^4
   e^{\frac{2}{3} \pi  k (5-3 q)}-3 e^{\frac{4}{3} \pi  k (3 q+2)}\right)}\\
   & =  \frac{4 \pi  \mu_a(H) (\langle \operatorname{w}^H, a\rangle-2 \pi  \mu_a(H) )}{3 \langle \operatorname{w}^H, a\rangle}.
\end{align*}
In particular, this is independent of $\xi$. Similarly, in the case of the exceptional divisor $E$, we have
\begin{align*}
&\frac{\vartheta_E }{(x y)^2\det \nabla^2 W }4\pi \mu_a(E) \big(\vartheta_{[\omega_k]} - \pi \frac{ \mu_a(E) }{\langle \operatorname{w}^E, a\rangle} W_k  - \pi \frac{ \mu_a(E) }{\langle \operatorname{w}^E, a\rangle} \vartheta_{E,k} \big)|_p \\
& = \frac{4 \pi  \mu_a(E)  x^4 \left(x^3 e^{2 \pi  k (1-3 q)} (q \langle \operatorname{w}^E, a\rangle x+2 \pi  \mu  (x+1))+\pi  \mu_a(E) -\langle \operatorname{w}^E, a\rangle\right)}{\langle \operatorname{w}^E, a\rangle \left(x^8 e^{2
   \pi  k (1-3 q)}-3 e^{2 \pi  k (3 q-1)}+2 x^4\right)},
\end{align*}
with leading contribution corresponding to $x \to -1$, given by
\begin{align*}
\lim_{k\to \infty}\frac{4 \pi  \mu_a(E)  \left(\langle \operatorname{w}^E, a\rangle \left(q e^{2 \pi  k (1-3 q)}-1\right)+\pi  \mu_a(E) \right)}{\langle \operatorname{w}^E, a\rangle \left(e^{2 \pi  k (1-3 q)}-3
   e^{2 \pi  k (3 q-1)}+2\right)} = 4 \pi  \mu_a(E) q.
\end{align*}
\addcontentsline{toc}{section}{References}
 
\bibliographystyle{abbrv}
 \bibliography{biblio_Kstab}

\begin{thebibliography}{10}

\bibitem{Arguz_equations}
H.~Arg\"{u}z.
\newblock The quantum mirror to the quartic del {P}ezzo surface.
\newblock In {\em Lie theory and its applications in physics}, volume 396 of
  {\em Springer Proc. Math. Stat.}, pages 423--429. Springer, Singapore, 2022.

\bibitem{Auroux_delPezzo}
D.~Auroux, L.~Katzarkov, and D.~Orlov.
\newblock Mirror symmetry for del {P}ezzo surfaces: vanishing cycles and
  coherent sheaves.
\newblock {\em Invent. Math.}, 166(3):537--582, 2006.

\bibitem{Barrott_equations}
L.~J. Barrott.
\newblock Explicit equations for mirror families to log {C}alabi-{Y}au
  surfaces.
\newblock {\em Bull. Korean Math. Soc.}, 57(1):139--165, 2020.

\bibitem{Chan_toric}
K.~Chan, S.-C. Lau, N.~C. Leung, and H.-H. Tseng.
\newblock Open {G}romov-{W}itten invariants and mirror maps for semi-{F}ano
  toric manifolds.
\newblock {\em Pure Appl. Math. Q.}, 16(3):675--720, 2020.

\bibitem{CheltsovRubinstein_flops}
I.~A. Cheltsov and Y.~A. Rubinstein.
\newblock On flops and canonical metrics.
\newblock {\em Ann. Sc. Norm. Super. Pisa Cl. Sci. (5)}, 18(1):283--311, 2018.

\bibitem{CoatesCortiIritani_hodge}
T.~Coates, A.~Corti, H.~Iritani, and H.-H. Tseng.
\newblock Hodge-theoretic mirror symmetry for toric stacks.
\newblock {\em J. Differential Geom.}, 114(1):41--115, 2020.

\bibitem{Dervan_crit_metrics}
R.~Dervan.
\newblock Stability conditions for polarised varieties.
\newblock {\em Forum Math. Sigma}, 11:Paper No. e104, 57, 2023.

\bibitem{DervanRoss_Kstab}
R.~Dervan and J.~Ross.
\newblock K-stability for {K}\"{a}hler manifolds.
\newblock {\em Math. Res. Lett.}, 24(3):689--739, 2017.

\bibitem{Donaldson_stabilitySurvey}
S.~K. Donaldson.
\newblock Stability of algebraic varieties and {K}\"{a}hler geometry.
\newblock In {\em Algebraic geometry: {S}alt {L}ake {C}ity 2015}, volume~97 of
  {\em Proc. Sympos. Pure Math.}, pages 199--221. Amer. Math. Soc., Providence,
  RI, 2018.

\bibitem{Givental_toric}
A.~Givental.
\newblock A mirror theorem for toric complete intersections.
\newblock In {\em Topological field theory, primitive forms and related topics
  ({K}yoto, 1996)}, volume 160 of {\em Progr. Math.}, pages 141--175.
  Birkh\"{a}user Boston, Boston, MA, 1998.

\bibitem{GriffithsHarris}
P.~Griffiths and J.~Harris.
\newblock {\em Principles of algebraic geometry}.
\newblock John Wiley \& Sons, 1978.

\bibitem{GrossHackingKeel_LCY}
M.~Gross, P.~Hacking, and S.~Keel.
\newblock Mirror symmetry for log {C}alabi-{Y}au surfaces {I}.
\newblock {\em Publ. Math. Inst. Hautes \'Etudes Sci.}, (122):65--168, 2015.

\bibitem{GrossHackingSiebert_theta}
M.~Gross, P.~Hacking, and B.~Siebert.
\newblock Theta functions on varieties with effective anti-canonical class.
\newblock {\em Mem. Amer. Math. Soc.}, 278(1367):xii+103, 2022.

\bibitem{Hertling}
C.~Hertling.
\newblock {\em Frobenius manifolds and moduli spaces for singularities}, volume
  151 of {\em Cambridge Tracts in Mathematics}.
\newblock Cambridge University Press, Cambridge, 2002.

\bibitem{KatzKontPant}
L.~Katzarkov, M.~Kontsevich, and T.~Pantev.
\newblock Bogomolov-{T}ian-{T}odorov theorems for {L}andau-{G}inzburg models.
\newblock {\em J. Differential Geom.}, 105(1):55--117, 2017.

\bibitem{Legendre_localised}
E.~Legendre.
\newblock Localizing the {D}onaldson-{F}utaki invariant.
\newblock {\em Internat. J. Math.}, 32(8):Paper No. 2150055, 23, 2021.

\bibitem{LiXu}
C.~Li and C.~Xu.
\newblock Special test configuration and {K}-stability of {F}ano varieties.
\newblock {\em Annals of Mathematics}, 180:197--232, 2014.

\bibitem{Odaka_blowup}
Y.~Odaka.
\newblock A generalization of the {R}oss-{T}homas slope theory.
\newblock {\em Osaka J. Math.}, 50(1):171--185, 2013.

\bibitem{RossThomas_obstruction}
J.~Ross and R.~Thomas.
\newblock An obstruction to the existence of constant scalar curvature
  {K}{\"a}hler metrics.
\newblock {\em Journal of Differential Geometry}, 72(3):429--466, 03 2006.

\bibitem{ScarpaStoppa_complexified}
C.~Scarpa and J.~Stoppa.
\newblock Special representatives of complexified {K}\"{a}hler classes.
\newblock {\em Selecta Math. (N.S.)}, 30(4):Paper No. 64, 45, 2024.

\bibitem{Wang_GIT}
X.~Wang.
\newblock Height and {GIT} weight.
\newblock {\em Math. Res. Lett.}, 19(4):909--926, 2012.

\bibitem{WangZhu_toric}
X.-J. Wang and X.~Zhu.
\newblock K\"{a}hler-{R}icci solitons on toric manifolds with positive first
  {C}hern class.
\newblock {\em Adv. Math.}, 188(1):87--103, 2004.

\end{thebibliography}

\noindent SISSA, via Bonomea 265, 34136 Trieste, Italy\\
Institute for Geometry and Physics (IGAP), via Beirut 2, 34151 Trieste, Italy\\
jstoppa@sissa.it    
\end{document}